\documentclass[twocolumn]{autart}

\usepackage{graphicx,algorithm,algorithmicx,algpseudocode,amsmath,amsfonts,verbatim,mathtools,enumerate,bm,hyperref,mathabx,hhline,multirow,tikz,braket,multirow,adjustbox,aircraftshapes,breakcites}
\usepackage[RPvoltages]{circuitikz}
\usepackage[style=base]{caption}
\usepackage{subcaption}
\usetikzlibrary{patterns}
\usetikzlibrary{calc}
\usepackage[flushleft]{threeparttable}

\allowdisplaybreaks

\newcommand{\argmin}{\mathop{\rm argmin}}

\newcommand{\norm}[1]{\left\lVert#1\right\rVert}
\newcommand{\mnorm}[1]{{\left\vert\kern-0.25ex\left\vert\kern-0.25ex\left\vert #1 
    \right\vert\kern-0.25ex\right\vert\kern-0.25ex\right\vert}}

\newtheorem{theorem}{Theorem}
\newtheorem{lemma}{Lemma}

\newtheorem{remark}{Remark}
\newtheorem{proposition}{Proposition}
\newtheorem{assumption}{Assumption}

\newcommand{\ie}{{\it i.e.}}

\begin{document}
\begin{frontmatter}
\runtitle{run title}

\title{Proportional-Integral Projected Gradient Method for Conic Optimization}

\author{Yue Yu},
\author{Purnanand Elango},
\author{Ufuk Topcu},
\author{Beh\c{c}et A\c{c}\i kme\c{s}e},
\address{Oden Institute for Computational Engineering and Sciences, The University of Texas, Austin, TX, 78712
}
\address{Department of Aeronautics and Astronautics, University of Washington, Seattle, WA, 98195
}

\thanks{Y. Yu and U. Topcu are with the Oden Institute for Computational Engineering and Sciences, The University of Texas, Austin, TX, 78712 (e-mails:  yueyu@utexas.edu,utopcu@utexas.edu). P. Elango and B. A\c{c}\i kme\c{s}e are with the William E. Boeing Department of Aeronautics \& Astronautics, University of Washington, Seattle, Washington, 98195 (e-mails:  pelango@uw.edu,behcet@uw.edu,
 }


\begin{abstract} 
Conic optimization is the minimization of a differentiable convex objective function subject to conic constraints. We propose a novel primal-dual first-order method for conic optimization, named proportional-integral projected gradient method (PIPG). PIPG ensures that both the primal-dual gap and the constraint violation converge to zero at the rate of \(O(1/k)\), where \(k\) is the number of iterations. If the objective function is strongly convex, PIPG improves the convergence rate of the primal-dual gap to \(O(1/k^2)\). Further, unlike any existing first-order methods, PIPG also improves the convergence rate of the constraint violation to \(O(1/k^3)\). We demonstrate the application of PIPG in constrained optimal control problems.
\end{abstract}

\end{frontmatter}

\section{Introduction}
Conic optimization is the minimization of a differentiable convex objective function subject to conic constraints:
\begin{equation}\label{opt: conic}
    \begin{array}{ll}
    \underset{z}{\mbox{minimize}} & f(z) \\
    \mbox{subject to} & Hz-g\in\mathbb{K}, \enskip z\in \mathbb{D},
    \end{array}
\end{equation}
where \(z\in\mathbb{R}^n\) is the solution variable,  \(f:\mathbb{R}^n\to\mathbb{R}\) is a continuously differentiable and convex objective function, \(\mathbb{K}\subset\mathbb{R}^m\) is a closed convex cone and \(\mathbb{D}\subset\mathbb{R}^n\) is a closed convex set, \(H\in\mathbb{R}^{m\times n}\) and \(g\in\mathbb{R}^m\) are constraint parameters. By proper choice of cone \(\mathbb{K}\), conic optimization \eqref{opt: conic} generalizes linear programming, quadratic programming, second-order cone programming, and semi-definite programming \cite{ben2001lectures,boyd2004convex}. 
Conic optimization has found applications in various areas, including signal processing \cite{luo2006introduction}, machine learning \cite{andersen2011interior}, robotics \cite{majumdar2020recent}, and aerospace engineering \cite{liu2017survey,eren2017model,malyuta2021advances}. 

The goal of numerically solving optimization \eqref{opt: conic} is to compute a solution \(z^\star\) that achieves, up to a given numerical tolerance, zero violation of the constraints in \eqref{opt: conic} and zero \emph{primal-dual gap}; the latter implies that \(z^\star\) is an optimal solution of optimization \eqref{opt: conic} \cite{boyd2011distributed,he20121,chambolle2011first,chambolle2016ergodic}. To this end, numerical methods iteratively compute a solution whose constraint violation and primal-dual gap are nonzero at first but converge to zero as the number of iteration \(k\) increases. 

Due to their low computational cost, first-order methods have attracted increasing attention in conic optimization \cite{lan2011primal,boyd2011distributed,o2016conic,chambolle2016introduction,yu2020proportional}. Unlike second-order methods, such as interior point methods \cite{nesterov1994interior,andersen2003implementing}, first-order methods do not rely on computing matrix inverses. They consequently are suitable for implementation  with limited computational resources.

\begin{table*}[!ht]
\centering
\begin{threeparttable}
\caption{Comparison of different first-order methods for conic optimization \eqref{opt: conic}} \label{tab: comparison}
  \begin{tabular}{|*{9}{c|}}
  \hline
  \multirow{3}{*}{Algorithms} & \multicolumn{4}{c|}{\(f\) is smooth \& convex} & \multicolumn{4}{c|}{\(f\) is smooth \& strongly convex}\\
  \cline{2-9}
  & \multicolumn{2}{c|}{\# of proj. per iter.} & \multicolumn{2}{c|}{convergence rates} & \multicolumn{2}{c|}{\# of proj. per iter.} & \multicolumn{2}{c|}{convergence rates}\\
  \cline{2-9}
  & \(\mathbb{D}\) & \(\mathbb{K}\) or \(\mathbb{K}^\circ\) & \shortstack{primal-dual \\ gap} & \shortstack{constraint \\ violation} & \(\mathbb{D}\) & \(\mathbb{K}\) or \(\mathbb{K}^\circ\) & \shortstack{primal-dual \\ gap} & \shortstack{constraint \\ violation}\\
  \hline
  ADMM & \(O(1/\sqrt{\epsilon})\) & 1 & \(O(1/k)\) & \(O(1/k)\) & \(O(\ln (1/\epsilon))\) & 1 & \(O(1/k)\) & \(O(1/k)\)\\
  PIPGeq & 1 & 1 & \(O(1/k)\) & \(O(1/k)\) & 1 & 1 & \(O(1/k)\) & \(O(1/k)\)\\
  PDHG & 1 & 1 & \(O(1/k)\) & N/A & 1 & 1 & \(O(1/k^2)\) & N/A\\
  This work & 1 & 1 & \(O(1/k)\) & \(O(1/k)\) & 1 & 1 & \(O(1/k^2)\) & \(O(1/k^3)\)\\
  \hline
  \end{tabular}
   \begin{tablenotes}
      \scriptsize
      \item \(\epsilon>0\) is a tunable accuracy tolerance in ADMM, \(\mathbb{K}^\circ\) denotes the polar cone of \(\mathbb{K}\). 
    \end{tablenotes}
  \end{threeparttable}
\end{table*}

The existing first-order methods solve optimization \eqref{opt: conic} by solving two different equivalent problems. The first equivalent problem is the following optimization with equality constraints \cite{boyd2011distributed,o2016conic,stellato2020osqp,yu2020proportional}:
\begin{equation}\label{opt: slack}
    \begin{array}{ll}
    \underset{z, y}{\mbox{minimize}} & f(z) \\
    \mbox{subject to} & Hz-y=g, \enskip y\in\mathbb{K}, \enskip z\in \mathbb{D}.
    \end{array}
\end{equation}
In particular, the alternating direction method of multipliers (ADMM) solves optimization \eqref{opt: conic} by computing one projection onto cone \(\mathbb{K}\) and multiple projections onto set \(\mathbb{D}\) in each iteration. ADMM ensures that both the constraint violation and the primal-dual gap converge to zero at rate of \(O(1/k)\), where \(k\) is the number of iterations \cite{gabay1976dual,eckstein1989splitting,fortin2000augmented,boyd2011distributed,he20121,wang2014bregman}. The proportional-integral projected gradient method for equality constrained optimization (PIPGeq) achieves the same convergence rates as ADMM, while computing one projection onto cone \(\mathbb{K}\) and \emph{only one} projection onto set \(\mathbb{D}\) in each iteration \cite{yu2020proportional}. Although variants of ADMM \cite{goldstein2014fast,kadkhodaie2015accelerated,ouyang2015accelerated,xu2017accelerated} and PIPGeq \cite{xu2017accelerated,yu2020proportional} can achieve accelerated convergence rates for strongly convex objective functions, such accelerations is not possible for optimization \eqref{opt: slack} because the objective function in \eqref{opt: slack} is independent of variable \(y\) and, as a result, not strongly convex. 

Another problem equivalent to optimization \eqref{opt: conic} is the following saddle-point problem, where \(\mathbb{K}^\circ\) is the polar cone of \(\mathbb{K}\) \cite{chambolle2011first,chambolle2016ergodic}:
\begin{equation}\label{opt: saddle}
    \begin{array}{lll}
        \underset{z\in\mathbb{D}}{\mbox{minimize}} & \underset{w\in\mathbb{K}^\circ}{\mbox{maximize}} & f(z)+\langle Hz-g, w\rangle.
    \end{array}
\end{equation}
In particular, the primal-dual hybrid-gradient method (PDHG) solves saddle-point problem \eqref{opt: saddle} by computing one projection onto cone \(\mathbb{K}^\circ\) and one projection onto set \(\mathbb{D}\) in each iteration. PDHG ensures that the primal-dual gap converges to zero at the rate of \(O(1/k)\) when for convex \(f\), and at an accelerated rate of \(O(1/k^2)\) for strongly convex \(f\) \cite{chambolle2016introduction,chambolle2016ergodic}.  However, since the constraint \(Hz-g\in\mathbb{K}\) is not explicitly considered in \eqref{opt: saddle}, the existing convergence results on PDHG do not provide any convergence rates of the violation of this constraint \cite{chambolle2016introduction,chambolle2016ergodic}.

We compare the per-iteration computation and the convergence rates of ADMM, PIPGeq and PDHG in Tab.~\ref{tab: comparison}. None of these methods \emph{simultaneously} has accelerated convergence rates (\ie, better than \(O(1/k)\)) for strongly convex \(f\) and guaranteed convergence rates on the constraint violation. To our best knowledge, whether there exists a first-order method that achieves both convergence results remains an open question.  

We answer this question affirmatively by proposing a novel primal-dual first-order method for conic optimization, named proportional-integral projected gradient method (PIPG). By combining the idea of proportional-integral feedback control and projected gradient method, PIPG ensures the following convergence results.
\begin{enumerate}
    \item For convex \(f\), both the primal-dual gap and the constraint violation converge to zero at the rate of \(O(1/k)\).
    \item For strongly convex \(f\), the convergence rate can be improved to \(O(1/k^2)\) for the primal-dual gap and \(O(1/k^3)\) for the constraint violation.
\end{enumerate}
PIPG generalizes both PDHG with constant step sizes \cite[Alg. 1]{chambolle2016ergodic} and PIPGeq \cite{yu2020proportional}. Compared with the existing methods, PIPG has the following advantages; see Tab.~\ref{tab: comparison} for an overview. In terms of per-iteration cost, it computes one projection onto cone \(\mathbb{K}^\circ\) and one projection onto set \(\mathbb{D}\), which is the same as PIPGeq and PDHG and fewer times of projections than ADMM. In terms of its convergence rates, to our best knowledge, the \(O(1/k^3)\) convergence rate of constraint violation has never been achieved before for general conic optimization. We numerically demonstrate these advantages of PIPG on several constrained optimal control problems.

The rest of the paper is organized as follows. After some preliminary results on convex analysis, Section~\ref{sec: preliminary} reviews existing first-order conic optimization methods. Section~\ref{sec: method} introduces PIPG along with its convergence results. Section~\ref{sec: numeric} demonstrate PIPG via numerical examples on constrained optimal control. Finally, Section~\ref{sec: conclusion} concludes and comments on future work.

\section{Preliminaries and related work}
\label{sec: preliminary}
This section reviews some basic results in convex analysis and several existing first-order conic optimization methods.

\subsection{Notation and preliminaries}\label{subsec: preliminary}

We let \(\mathbb{N}\), \(\mathbb{R}\) and \(\mathbb{R}_+\) denote the set of positive integer, real, and non-negative real numbers, respectively. For two vectors \(z, z'\in\mathbb{R}^n\), \(\langle z, z'\rangle\) denotes their inner product, \(\norm{z}\coloneqq \sqrt{\langle z, z\rangle}\) denotes the \(\ell_2\) norm of \(z\), and \(\norm{\cdot}_\infty\) denotes the \(\ell_\infty\) norm of \(z\), \ie, the maximum absolute value of the entries of \(z\). We let \(\mathbf{1}_n\) and \(\mathbf{0}_n\) denote the \(n\)-dimensional vectors of all \(1\)'s and all \(0\)'s, respectively. We also let \(I_n\) and \(0_{m\times n}\) denote the \(n\times n\) identity matrix and the \(m\times n\) zero matrix, respectively. When their dimensions are clear from the context, we omit the subscripts and simply write vector \(\mathbf{1}, \mathbf{0}\) and matrix \(I,0\). For a matrix \(H\in\mathbb{R}^{m\times n}\), \(H^\top\) denotes its transpose, \(\mnorm{H}\) denotes its largest singular value. For a square matrix \(M\in\mathbb{R}^{n\times n}\), \(\exp(M)\) denotes the matrix exponential of \(M\), and \(\norm{z}_M\coloneqq \sqrt{\langle z, M z\rangle}\) for all  \(z\in\mathbb{R}^n\). Given two sets \(\mathbb{S}_1\) and \(\mathbb{S}_2\), \(\mathbb{S}_1\times\mathbb{S}_2\) denotes their Cartesian product.

Let \(z, z'\in\mathbb{R}^n\) and \(f:\mathbb{R}^n\to\mathbb{R}\) be a continuously differentiable function. The \emph{Bregman divergence} from \(z\) to \(z'\) associated with function \(f\) is given by
\begin{equation}\label{eqn: Bregman}
    B_f(z, z')\coloneqq f(z)-f(z')-\langle\nabla f(z'), z-z'\rangle.
\end{equation}
We say function \(f\) is \(\mu\)-strongly convex for some \(\mu\in\mathbb{R}_+\) if
\begin{equation}\label{eqn: strongly convex}
    B_f(z, z')\geq \frac{\mu}{2}\norm{z-z'}^2
\end{equation}
for all \(z, z'\in\mathbb{R}^n\). When \eqref{eqn: strongly convex} holds with \(\mu=0\), we say function \(f\) is convex. We say function \(f\) is \(\lambda\)-smooth  for some \(\lambda\in\mathbb{R}_+\) if
\begin{equation}
    B_f(z, z')\leq \frac{\lambda}{2}\norm{z-z'}^2
\end{equation}
for all \(z, z'\in\mathbb{R}^n\).

Let \(\mathbb{D}\subset\mathbb{R}^n\) be a closed convex set, \ie, \(\mathbb{D}\) contains all of its boundary points and \(\gamma z+(1-\gamma)z'\in\mathbb{D}\) for any \(\gamma\in[0, 1]\) and \(z, z'\in\mathbb{D}\). The \emph{projection} of \(z\in\mathbb{R}^n\) onto set \(\mathbb{D}\) is given by
\begin{equation}\label{eqn: projection}
\pi_{\mathbb{D}}[z]\coloneqq \underset{y\in\mathbb{D}}{\argmin} \norm{z-y}.
\end{equation}
Let \(\mathbb{K}\subset\mathbb{R}^m\) be a closed convex cone, \ie, \(\mathbb{K}\) is a closed convex set and \(\gamma w\in\mathbb{K}\) for any \(w\in\mathbb{K}\) and \(\gamma\in\mathbb{R}_+\). The \emph{polar cone} of \(\mathbb{K}\) is also a closed convex cone given by
\begin{equation}\label{eqn: polar cone}
    \mathbb{K}^\circ\coloneqq \{w\in\mathbb{R}^m| \langle w, y\rangle\leq 0, \forall y\in\mathbb{K}\}. 
\end{equation}

\subsection{Related work}
\label{subsec: related work}
We briefly review three existing first-order primal-dual conic optimization methods: ADMM, PIPGeq, and PDHG. In the following, let \(\alpha, \beta, \gamma\) denote positive scalar step sizes, and \(\{\alpha^j\}_{j\in\mathbb{N}}\), \(\{\beta^j\}_{j\in\mathbb{N}}\), \(\{\gamma^j\}_{j\in\mathbb{N}}\) denote sequences of positive scalar step sizes. For simplicity, we assume all methods are terminated after a fixed number of iterations, denoted by \(k\in\mathbb{N}\).

\subsubsection{Alternating direction method of multipliers}
As a special case of Douglas-Rachford splitting method \cite{eckstein1989splitting,fortin2000augmented}, alternating direction method of multipliers (ADMM) solves optimization \eqref{opt: conic} by solving the equivalent optimization \eqref{opt: slack} using Algorithm~\ref{alg: ADMM} \cite{gabay1976dual,boyd2011distributed,he20121}.
\begin{algorithm}[!ht]
\caption{ADMM}
\begin{algorithmic}[1]
\Require \(k, \alpha, z^1\in\mathbb{D}, y^1\in\mathbb{K}, w^1\in\mathbb{R}^m\)
\Ensure \(z^{k}\)
\For{\(j=1, 2, \ldots, k-1\)}
\State{\(z^{j+1}=\underset{z\in\mathbb{D}}{\argmin}\, f(z)+\frac{\alpha}{2}\norm{Hz-y^j-g+w^j}^2\)}\label{eqn: ADMM prox}
\State{\(y^{j+1}=\pi_{\mathbb{K}}[Hz^{j+1}-g+w^j]\)}
\State{\(w^{j+1}=w^j+Hz^{j+1}-y^{j+1}-g\)}
\EndFor
\end{algorithmic}
\label{alg: ADMM}
\end{algorithm}

Generally, the minimization in the line ~\ref{eqn: ADMM prox} of Algorithm~\ref{alg: ADMM} can only be solved approximately up to a numerical tolerance \(\epsilon>0\) using iterative methods. Such methods need to compute at least \(O(1/\sqrt{\epsilon})\) projections onto set \(\mathbb{D}\) if \(f\) is merely convex, and \(O(\ln \frac{1}{\epsilon})\) projections if function \(f\) is strongly convex; see \cite[Chp. 2]{nesterov2018lectures} for a detailed discussion.

There has been many variants of ADMM developed in the literature. However, none of them lead to any significant benefits for optimization in \eqref{opt: slack}. For example, \cite{ouyang2015accelerated} and \cite[Alg. 1]{xu2017accelerated} simplified the minimization in the line \eqref{eqn: ADMM prox} of Algorithm~\ref{alg: ADMM} by approximating function \(f\) using its linearization. However, solving the resulting approximate minimization still  requires multiple projections onto set \(\mathbb{D}\). On the other hand, although the convergence of ADMM can be accelerated when the objective function is strongly convex \cite{goldstein2014fast,kadkhodaie2015accelerated,ouyang2015accelerated,xu2017accelerated}, such acceleration does not apply to the optimization \eqref{opt: slack}. The reason is because the objective function in \eqref{opt: slack} is not strongly convex with respect to (in fact, does not depend on) variable \(y\). 

\subsubsection{Proportional-integral projected gradient method for equality constrained optimization}
Motivated by applications in model predictive control, the proportional-integral projected gradient method for equality constrained optimization (PIPGeq) solves optimization \eqref{opt: conic} by solving the equivalent optimization \eqref{opt: slack} using Algorithm~\ref{alg: PIPGeq}.

\begin{algorithm}[!ht]
\caption{PIPGeq}
\begin{algorithmic}[1]
\Require \(k, \alpha,\beta, z^1\in\mathbb{D}, y^1\in\mathbb{K}, w^1\in\mathbb{R}^m\).
\Ensure \(z^{k}\).
\For{\(j=1, 2, \ldots, k-1\)}
\State{\(v^{j+1}=w^j+\beta(Hz^j-y^j-g)\)}
\State{\(z^{j+1}=\pi_{\mathbb{D}}[z^j-\alpha(\nabla f(z^j)+H^\top v^{j+1})]\)}\label{eqn: PIPGeq proj}
\State{\(y^{j+1}=\pi_{\mathbb{K}}[y^j+\alpha v^{j+1}]\)}
\State{\(w^{j+1}=w^j+\beta(Hz^{j+1}-y^{j+1}-g)\)}
\EndFor
\end{algorithmic}
\label{alg: PIPGeq}
\end{algorithm}

Unlike line~\ref{eqn: ADMM prox} in Algorithm~\ref{alg: ADMM}, line \ref{eqn: PIPGeq proj} in Algorithm~\ref{alg: PIPGeq} computes only one projection onto set \(\mathbb{D}\) instead of multiple times. As a result, PIPGeq can achieve the same convergence rates as those of ADMM while lowering the per-iteration computation cost \cite{xu2017accelerated,yu2020proportional}.

\subsubsection{Primal-dual hybrid gradient method}
Motivated by applications in computational imaging, the primal-dual hybrid gradient method (PDHG) was first introduced in \cite{chambolle2011first} and later shown to be equivalent to Douglas-Rachford splitting method \cite{o2020equivalence}. Later, another variant of PDHG was introduced in \cite{chambolle2016ergodic}, which is an instance of three-operator splitting methods \cite{vu2013splitting,condat2013primal,chen2016primal,davis2017three,yan2018new}. To solve optimization \eqref{opt: conic}, PDHG solves the equivalent convex-concave saddle point problem \eqref{opt: saddle} instead. If function \(f\) is merely convex, PDHG uses Algorithm~\ref{alg: PDHGc}. If function \(f\) is \(\mu\)-strongly convex for some \(\mu>0\), then PDHG uses Algorithm~\ref{alg: PDHGv} instead.

\begin{algorithm}[!ht]
\caption{PDHG with constant step sizes}
\begin{algorithmic}[1]
\Require \(k, \alpha, \beta, z^1\in\mathbb{D}, w^1\in\mathbb{K}^\circ\).
\Ensure \(z^{k}\).
\For{\(j=1, 2, \ldots, k-1\)}
\State{\(z^{j+1}=\pi_{\mathbb{D}}[z^j-\alpha(\nabla f(z^j)+H^\top w^j)]\)}
\State{\(w^{j+1}=\pi_{\mathbb{K}^\circ}[w^j+\beta (H(2z^{j+1}-z^j)-g)]\)}
\EndFor
\end{algorithmic}
\label{alg: PDHGc}
\end{algorithm}

\begin{algorithm}[!ht]
\caption{PDHG with varying step sizes}
\begin{algorithmic}[1]
\Require \(k, \{\alpha^j, \beta^j, \gamma^j\}_{j=1}^k, \mu, z^1\in\mathbb{D}, w^1\in\mathbb{K}^\circ\).
\Ensure \(z^{k}\).
\For{\(j=1, 2, \ldots, k-1\)}
\State{\(w^{j+1}=\pi_{\mathbb{K}^\circ}[w^j+\beta^j (H(z^j+\gamma^j(z^j-z^{j-1}))-g)]\)}
\State{\(z^{j+1}=\pi_{\mathbb{D}}\left[z^j-\frac{\alpha^j}{\mu\alpha^j+1}(\nabla f(z^j)+H^\top w^{j+1})\right]\)}
\EndFor
\end{algorithmic}
\label{alg: PDHGv}
\end{algorithm}

The primal-dual gap converges to zero at the rate of \(O(1/k)\) and \(O(1/k^2)\) for Algorithm~\ref{alg: PDHGc} and Algorithm~\ref{alg: PDHGv}, respectively \cite{chambolle2016ergodic}. However, to our best knowledge, there is no convergence result on the constraint violation for either Algorithm~\ref{alg: PDHGc} or Algorithm~\ref{alg: PDHGv}.
\section{Proportional-integral projected gradient method}
\label{sec: method}
We introduce a novel first-order primal-dual method, named \emph{proportional-integral projected gradient method (PIPG)}, for conic optimization \eqref{opt: conic}, and discuss its convergence rates in terms of the constraint violation and the primal-dual gap.

 Algorithm~\ref{alg: PIPG} summarizes the proposed method, where \(k\in\mathbb{N}\) is the maximum number of iterations, and \(\{\alpha^j\}_{j=1}^k\) and \(\{\beta^j\}_{j=1}^k\) are sequences of positive scalar step sizes that will be specified later. We note that, instead of maximum number of iterations, one can use alternative stopping criterions, such as the distance between \(Hz^j-g\) and \(\mathbb{K}\) reaching a given tolerance.
\begin{algorithm}[!ht]
\caption{PIPG}
\begin{algorithmic}[1]
\Require \(k, \{\alpha^j, \beta^j\}_{j=1}^k, z^1\in\mathbb{D}, v^1\in\mathbb{K}^\circ\).
\Ensure \(z^k\).
\For{\(j=1, 2, \ldots, k-1\)}
\State{\(w^{j+1}= \pi_{\mathbb{K}^\circ}[v^j+\beta^j(Hz^j-g)]\)}\label{alg: w}
\State{\(z^{j+1}=\pi_{\mathbb{D}}[z^j-\alpha^j(\nabla f(z^j)+H^\top w^{j+1})]\)}\label{alg: z}
\State{\(v^{j+1}=w^{j+1}+\beta^j H(z^{j+1}-z^j)\)}\label{alg: v}
\EndFor
\end{algorithmic}
\label{alg: PIPG}
\end{algorithm}

The name PIPG is due to the following observations. First, if \(\mathbb{K}=\{\mathbf{0}_n\}\),  then \(\mathbb{K}^\circ=\mathbb{R}^m\) and line~2 and line~4 in Algorithm~\ref{alg: PIPG} become the following:
\begin{subequations}\label{alg: PIPGeq 2}
\begin{align}
    w^{j+1}=& v^j+\beta^j(Hz^j-g),\label{alg: pro-integ} \\
    v^{j+1}=&v^j+\beta^j (Hz^{j+1}-g).\label{alg: integral}
\end{align}
\end{subequations}
Using \eqref{alg: integral} one can show that
\[v^j=v^1+\sum_{i=2}^j\beta^{i-1} (Hz^i-g), \]
Hence \(v^k\) is a weighted summation, or numerical \emph{integration}, of \(Hz^i-g\) from \(i=2\) to \(i=j\). Further, \eqref{alg: pro-integ} states that \(w^j\) adds a \emph{proportional} term of \(Hz^j-g\) to \(v^j\), hence \(w^j\) in \eqref{alg: pro-integ} is a \emph{proportional-integral} term of \(Hz^j-g\). Second, if \(H\) is a zero matrix, then line~\ref{alg: z} in Algorithm~\ref{alg: PIPG} becomes a projected gradient method that minimizes \(f\) over set \(\mathbb{D}\) \cite[Sec. 2.2.5]{nesterov2018lectures}. Therefore Algorithm~\ref{alg: PIPG} can be interpreted as a combination of proportional-integral feedback control and the projected gradient method. Similar idea has also been popular in equality constrained optimization \cite{wang2010control,yu2020mass,yu2020rlc,yu2020proportional}. 

\begin{remark}
Notice that the \(w^{j+1}\) in \eqref{alg: pro-integ} is otherwise identical to the \(v^{j+1}\) in \eqref{alg: integral} except that \eqref{alg: pro-integ} uses \(z^j\) whereas \eqref{alg: integral} uses \(z^{j+1}\). Such scheme is also known as a \emph{prediction-correction} step, which has been popular in many first-order primal-dual methods, including the extra-gradient and mirror-prox method \cite{korpelevich1977extragradient,nemirovski2004prox,nesterov2007dual}, the accelerated linearized ADMM \cite{ouyang2015accelerated,xu2017accelerated}, the primal-dual fixed point methods \cite{krol2012preconditioned,chen2013primal,chen2016primal,yan2018new} and the accelerated mirror descent method \cite{cohen2018acceleration}.
\end{remark}

\begin{remark}
One can verify that if \(\alpha^j\equiv \alpha\) and \(\beta^j\equiv\beta\) for \(j=1, 2, \ldots, k\), then Algorithm~\ref{alg: PIPG} is equivalent to Algorithm~\ref{alg: PDHGc}, the latter was first introduced in \cite[Alg. 1]{chambolle2016ergodic}. 
\end{remark}

Next, we will show the convergence results of Algorithm~\ref{alg: PIPG}. To this end, we will frequently use the following  \emph{quadratic distance function} to closed convex cone \(\mathbb{K}\):
\begin{equation}\label{eqn: dist func}
    d_{\mathbb{K}}(w)\coloneqq\underset{v\in\mathbb{K}}{\mbox{minimize}}\,\,\frac{1}{2}\norm{w-v}^2,
\end{equation}
which is continuously differentiable and convex \cite[Lem. 2.2.9]{nesterov2018lectures}. We will also use the following \emph{Lagrangian function}:
\begin{equation}\label{eqn: Lagrangian}
    L(z, w)\coloneqq f(z)+\langle Hz-g, w\rangle.
\end{equation}

We make the following assumptions on optimization \eqref{opt: conic}.

\begin{assumption}\label{asp: opt}
\begin{enumerate}
    \item Function \(f:\mathbb{R}^n\to\mathbb{R}\) is continuously differentiable. There exists \(\mu, \lambda\in\mathbb{R}_+\) with \(\mu\leq \lambda\) such that \(f\) is \(\mu\)-strongly convex and \(\lambda\)-smooth, \ie, 
    \begin{equation*}
        \frac{\mu}{2}\norm{z-z'}^2\leq B_f(z, z')\leq \frac{\lambda}{2}\norm{z-z'}^2
    \end{equation*}
    for all \(z, z'\in\mathbb{R}^n\).
    \item Set \(\mathbb{D}\subset\mathbb{R}^n\) and cone \(\mathbb{K}\subset\mathbb{R}^m\) are closed and convex. 
    \item There exists \(z^\star\in\mathbb{D}\) and \(w^\star\in\mathbb{K}^\circ\) such that
    \begin{equation*}
        L(z^\star, \overline{w})\leq L(z^\star, w^\star)\leq L(\overline{z}, w^\star)
    \end{equation*}
    for all \(\overline{z}\in\mathbb{D}\) and  \(\overline{w}\in\mathbb{K}^\circ\).
\end{enumerate}
\end{assumption}

Under the above assumptions, the quantity
\(L(\overline{z}, w^\star)-L(z^\star, \overline{w})\), also known as the \emph{primal-dual gap} evaluated at \((\overline{z},\overline{w})\), is non-negative \cite{boyd2011distributed,he20121,chambolle2011first,chambolle2016ergodic}. The following proposition provides a sufficient condition on \(\overline{z}\) and \(\overline{w}\) under which the primal-dual gap \(L(\overline{z}, w^\star)-L(z^\star, \overline{w})\) equals zero and \(\overline{z}\) is an optimal solution of optimization \eqref{opt: conic}. 

\begin{proposition}\label{prop: kkt}
If there exists \(\overline{z}\in\mathbb{D}\) and \(\overline{w}\in\mathbb{K}^\circ\) such that
\begin{equation}\label{eqn: suboptimality}
    L(\overline{z}, w)-L(z, \overline{w})\leq 0,
\end{equation}
for all \(z\in\mathbb{D}\) and \(w\in\mathbb{K}^\circ\), then \(\overline{z}\) is an optimal solution of optimization \eqref{opt: conic}, \ie,  \(H\overline{z}-g\in\mathbb{K}\) and \(f(\overline{z})\leq f(z)\) for any \(z\in\mathbb{D}\) such that \(Hz-g\in\mathbb{K}\).
\end{proposition}

\begin{pf*}{Proof}
See Appendix~\ref{app: prop}.
\end{pf*}

As our first step, the following lemma proves a key inequality for our later discussions.

\begin{lemma}\label{lem: iteration}
Suppose that Assumption~\ref{asp: opt} holds and \(\{w^j, z^j, v^j\}_{j=1}^k\) is computed using Algorithm~\ref{alg: PIPG} where \(\alpha^j, \beta^j>0\) and \(\alpha^j(\lambda +\sigma \beta^j)=1\) for some \(\sigma\geq \mnorm{H}^2\) and all \(j=1, 2, \ldots k\)  . Then
\begin{equation*}
\begin{aligned}
    &\beta^jd_{\mathbb{K}}(Hz^j-g)+L(z^{j+1}, w)-L(z, w^{j+1})\\
    \leq & \left(\frac{1}{2\alpha^j}-\frac{\mu}{2}\right)\norm{z^j-z}^2+ \frac{1}{2\beta^j}\norm{v^j-w}^2\\
    &-\frac{1}{2\alpha^j}\norm{z^{j+1}-z}^2-\frac{1}{2\beta^j}\norm{v^{j+1}-w}^2,
\end{aligned}    
\end{equation*}
for all \(z\in\mathbb{D}\), \(w\in\mathbb{K}^\circ\), and \(j=1, 2, \ldots, k\).
\end{lemma}
\begin{pf*}{Proof}
See Appendix~\ref{app: lem}.
\end{pf*}

Equipped with Lemma~\ref{lem: iteration}, we are ready to prove the convergence results of Algorithm~\ref{alg: PIPG}. The idea is to first summing up the inequality in Lemma~\ref{lem: iteration} corresponding to different value of \(j\), then using the Jensen's inequality.

We start with the case where \(\mu=0\), \ie, function \(f\) is merely convex. The following theorem shows the convergence results of Algorithm~\ref{alg: PIPG} in this case.

\begin{theorem}\label{thm: convex}
Suppose that Assumption~\ref{asp: opt} holds with \(\mu=0\), and \(\{w^j, z^j, v^j\}_{j=1}^k\) is computed using Algorithm~\ref{alg: PIPG} with \(\alpha^j=\frac{1}{\beta\sigma+\lambda}\) and \(\beta^j=\beta\) and all \(j=1, 2, \ldots, k\), where \(\beta>0\) and \(\sigma\geq \mnorm{H}^2\). Let 
\begin{equation*}
    \tilde{z}^k\coloneqq\frac{1}{k}\sum_{j=1}^k z^j,\enskip \overline{z}^k\coloneqq \frac{1}{k}\sum_{j=1}^k z^{j+1},\enskip \overline{w}^k\coloneqq\frac{1}{k}\sum_{j=1}^k w^{j+1}, 
\end{equation*}
and \(V^1(z, w)\coloneqq\frac{1}{2\alpha }\norm{z^1-z}^2+\frac{1}{2\beta}\norm{v^1-w}^2\) for all \(z\in\mathbb{D}\) and \(w\in \mathbb{K}^\circ\). Then \(\tilde{z}^k, \overline{z}^k\in\mathbb{D}\), \(\overline{w}^k\in\mathbb{K}^\circ\), and
\begin{equation*}
\begin{aligned}
      d_{\mathbb{K}}(H \tilde{z}^k-g)
    \leq & \frac{V^1(z^\star, w^\star)}{\beta k}, \\
        L( \overline{z}^k,w)-L(z,  \overline{w}^k)
    \leq  & \frac{V^1(z, w)}{k},
\end{aligned}
\end{equation*}
for all \(z\in\mathbb{D}\), \(w\in\mathbb{K}^\circ\).
\end{theorem}

\begin{pf*}{Proof}
See Appendix~\ref{app: thm 1}.
\end{pf*}

Theorem~\ref{thm: convex} shows that \(\tilde{z}^k,\overline{z}^k\in\mathbb{D}\). In addition, as \(k\) increases, the violation of constraint \(H\tilde{z}^k-g\in\mathbb{K}\), measured by nonnegative distance \(d_{\mathbb{K}}(H\tilde{z}^k-g)\), converges to zero, and the condition in \eqref{eqn: suboptimality} holds asymptotically for \(\overline{z}^k\) and \(\overline{w}^k\).

If \(\mu>0\), \ie, function \(f\) is strongly convex, then  we can further improve the convergence results in Theorem~\ref{thm: convex} as follows.

\begin{theorem}\label{thm: strongly convex}
Suppose that Assumption~\ref{asp: opt} holds with \(\mu>0\) and
\(\{w^j, z^j, v^j\}_{j=1}^k\) is computed using Algorithm~\ref{alg: PIPG} with \(\alpha^j=\frac{2}{(j+1)\mu+2\lambda}\) and \(\beta^j=\frac{(j+1)\mu}{2\sigma}\) for some \(\sigma>\mnorm{H}^2\) and all \(j=1,2,\ldots, k\). Let
\begin{equation*}
\begin{aligned}
       \tilde{z}^k\coloneqq & \frac{3}{k(k^2+6k+11)}\sum_{j=1}^k (j+1)(j+2) z^j,\\
      \overline{z}^k\coloneqq& \frac{2}{k(k+5)}\sum_{j=1}^k (j+2)z^{j+1},\\ \overline{w}^k\coloneqq& \frac{2}{k(k+5)}\sum_{j=1}^k (j+2)w^{j+1},
\end{aligned}
\end{equation*}
and \(V^1(z, w)\coloneqq\frac{\mu+2\lambda}{4 }\norm{z^1-z}^2+\frac{\sigma}{\mu}\norm{v^1-w}^2\) for all \(z\in\mathbb{D}\) and \(w\in \mathbb{K}^\circ\). Then \(\tilde{z}^k, \overline{z}^k\in\mathbb{D}\), \(\overline{w}^k\in\mathbb{K}^\circ\), and
\begin{equation*}
\begin{aligned}
      d_{\mathbb{K}}(H \tilde{z}^k-g)\leq &  \frac{12\lambda\sigma V^1(z^\star, w^\star)}{\mu^2 k(k^2+6k+11)},\\
      L( \overline{z}^k,w)-L(z,  \overline{w}^k)
    \leq & \frac{4\lambda V^1(z, w)}{\mu k(k+5)},
\end{aligned}
\end{equation*}
for all \(z\in\mathbb{D}\), \(w\in\mathbb{K}^\circ\).
\end{theorem}

\begin{pf*}{Proof}
See Appendix~\ref{app: thm 2}.
\end{pf*}

\begin{remark}
Unlike the results in \cite{chambolle2016ergodic}, Theorem~\ref{thm: convex} and Theorem~\ref{thm: strongly convex} prove not only the convergence of the primal-dual gap, but also the convergence of the constraints violation. In addition, if \(\alpha^j\equiv \alpha\) and \(\beta^j\equiv\beta\) for \(j=1, 2, \ldots, k\), then one can show that Algorithm~\ref{alg: PIPG} is equivalent to Algorithm~\ref{alg: PDHGc}; in other words, the results in Theorem~\ref{thm: convex} also apply to Algorithm~\ref{alg: PDHGc}.  
\end{remark}

\begin{remark}
When using varying step sizes, Algorithm~\ref{alg: PIPG} differs from Algorithm~\ref{alg: PDHGv} in the relation between step sizes and the iteration number: the one in Algorithm~\ref{alg: PIPG} is explicit, whereas the one in Algorithm~\ref{alg: PDHGv} is implicitly defined by a recursive formula \cite[Sec. 5.2]{chambolle2016ergodic}. Furthermore, we can prove the convergence rate of the constraint violation for Algorithm~\ref{alg: PIPG}, whereas similar rate for Algorithm~\ref{alg: PDHGv} is, to our best knowledge, does not exist in the literature.
\end{remark}

\section{Applications to constrained optimal control}
\label{sec: numeric}
We demonstrate the application of PIPG to constrained optimal control problems. In Section~\ref{subsec: optimal control}, we show how to formulate a typical constrained optimal control problem as an instance of conic optimization \eqref{opt: conic}, and provide examples from mechanical engineering and robotics. In Section~\ref{subsec: numerics}, we demonstrate the performance of PIPG via said examples, and compare it against the existing methods reviewed in Section~\ref{subsec: related work}. Throughout we let \(n_x, n_u, p_x, p_u\in\mathbb{N}\) denote the dimension of different vector spaces, \(\Delta\in\mathbb{R}_+\) denote a positive sampling time period, and \(t\in\mathbb{N}\) denote a discrete time index.

\subsection{Constrained optimal control}\label{subsec: optimal control}

We consider the following linear time invariant system
\begin{equation}\label{sys: CT ATI}
    \frac{d}{ds}x(s)=A_cx(s)+B_cu(s)+h_c
\end{equation}
where \(x:\mathbb{R}_+\to\mathbb{R}^{n_x}\) and \(u: \mathbb{R}_+\to\mathbb{R}^{n_u}\) denote the state and input function, respectively. Matrix \(A_c\in\mathbb{R}^{n_x\times n_x}\), \(B_c\in\mathbb{R}^{n_x\times n_u}\), and vector \(h_c\in\mathbb{R}^{n_x}\) are known parameters. 

If the input changes value only at discrete time instants, then we can simplify dynamics \eqref{sys: CT ATI} as follows. Let \(\Delta\in\mathbb{N}_+\) and
\(x_t\coloneqq x(t\Delta), \enskip u_t\coloneqq u(t\Delta)\)
for all \(t\in\mathbb{N}\). Suppose that
\[u(s)=u(t\Delta), \enskip t\Delta\leq s<(t+1)\Delta,\]
for all \(t\in\mathbb{N}\). Then dynamics equation \eqref{sys: CT ATI} is equivalent to the following
\begin{equation}\label{sys: DT ATI}
    x_{t+1}=Ax_t+Bu_t+h,
\end{equation}
for all \(t\in\mathbb{N}\), where
\begin{equation}\label{eqn: ZOH}
\begin{aligned}
     &\textstyle A=\exp(A_c\Delta), \enskip B=\left(\int_0^\Delta \exp(A_c s)ds\right)B_c,\\
     & \textstyle h=\left(\int_0^\Delta \exp(A_c s)ds\right)h_c.
\end{aligned}
\end{equation}
For further details on the above equivalence, we refer the interested readers to \cite[Sec. 4.2.1]{chen1999linear}.

Let \(\{x_{t+1}, u_t\}_{t=0}^{\tau-1}\) denote a length-\(\tau\) input-state trajectory of system \eqref{sys: DT ATI} for some \(\tau\in\mathbb{N}\), and \(\{\hat{x}_{t+1},\hat{u}_t\}_{t=1}^{\tau-1}\) denote a desired length-\(\tau\) reference input-state trajectory. A typical optimal control problem is the minimization of the difference between \(\{x_{t+1},u_t\}_{t=0}^{\tau-1}\) and \(\{ \hat{x}_{t+1},\hat{u}_t\}_{t=0}^{\tau-1}\) subject to various constraints:
\begin{subequations}\label{opt: mpc}
    \begin{align}
        \underset{\{u_t, x_{t+1}\}_{t=0}^{\tau-1}}{\mbox{minimize}} \enskip & \textstyle \frac{1}{2}\sum\limits_{t=0}^{\tau-1}(\norm{x_{t+1}-\hat{x}_{t+1}}_Q^2+\norm{u_t-\hat{u}_t}_R^2) \label{mpc: cost}\\
        \mbox{subject to} \enskip &  x_{t+1}=Ax_t+Bu_t+h,\enskip 0\leq t\leq \tau-1,\label{mpc: dynamics}\\
        &\norm{u_{t+1}-u_t}_\infty\leq \gamma , \enskip 0\leq t\leq \tau-2,\label{mpc: input rate}\\
        & C_tx_t-a_t\geq 0,\, x_t\in\mathbb{X},\enskip 1\leq t\leq \tau,\label{mpc: state}\\
        & D_tu_t-b_t\geq 0,\, u_t\in\mathbb{U},\enskip 0\leq t\leq \tau-1.\label{mpc: input}
    \end{align}
\end{subequations}
In particular, the objective function in \eqref{mpc: cost} is a quadratic distance between \(\{x_{t+1},u_t\}_{t=0}^{\tau-1}\) and \(\{ \hat{x}_{t+1},\hat{u}_t\}_{t=0}^{\tau-1}\), where \(Q\in\mathbb{R}^{n_x\times n_x}\) and \(R\in\mathbb{R}^{n_u\times n_u}\) are given symmetric and positive definite weighting matrices. The constraints in \eqref{mpc: dynamics} ensure that inpiut-state trajectory \(\{u_{0:\tau-1}, x_{1:\tau}\}\) agree with the dynamics \eqref{sys: DT ATI}, where \(x_0\in\mathbb{R}^{n_x}\) is the given initial state.  The constraints in \eqref{mpc: input rate} upper bound the elementwise difference between two consecutive inputs by \(\gamma\in\mathbb{R}_+\), which prevents frequent and large input variations \cite[Sec. 4.10]{betts2010practical}. The constraints in \eqref{mpc: state} and \eqref{mpc: input} describe possible physical and operational constraints on states and inputs, where \(C_t\in\mathbb{R}^{p_x\times n_x}\), \(D_t\in\mathbb{R}^{p_u\times n_u}\), \(a_t\in\mathbb{R}^{p_x}\), \(b_t\in\mathbb{R}^{p_u}\), \(\mathbb{X}\subset\mathbb{R}^{n_x}\) and \(\mathbb{U}\subset\mathbb{R}^{n_u}\) are closed convex sets. 

One can transform optimization \eqref{opt: mpc} into a special case of optimization \eqref{opt: conic} using particular choices of the parameters. See Appendix~\ref{app: matrices} for the detailed transformation.

In the following, we will provide two illustrating examples of optimization \eqref{opt: mpc} from mechanical engineering and robotics applications. For simplicity, all problem parameters will be unitless.

\subsubsection{Oscillating masses control}
\begin{figure}[!ht]
\centering

\begin{adjustbox}{scale=0.6}

\begin{circuitikz}

\pattern[pattern=north east lines] (-0.2, -0.8) rectangle (0, 0.8);
\draw[thick] (0, -0.8) -- (0, 0.8);
\draw (0, 0) to[spring] (1, 0);
\draw (1, -0.5) rectangle (2, 0.5);
\draw (2, 0) to[spring] (3, 0);
\draw (3, -0.5) rectangle (4, 0.5);
\draw (4, 0) to[spring] (5, 0);

\node at (5.5, 0) {\huge $\ldots$};

\draw (6, 0) to[spring] (7, 0);
\draw (7, -0.5) rectangle (8, 0.5);
\draw (8, 0) to[spring] (9, 0);
\draw (9, -0.5) rectangle (10, 0.5);
\draw (10, 0) to[spring] (11, 0);

\pattern[pattern=north east lines] (11, -0.8) rectangle (11.2, 0.8);
\draw[thick] (11, -0.8) -- (11, 0.8);
\end{circuitikz}

\end{adjustbox}
\caption{The oscillating masses system}
\label{fig: mass-spring model}
\end{figure}
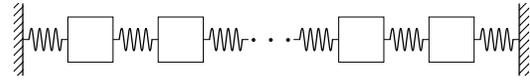

We consider the problem of controlling a one-dimensional oscillating masses system using external forcing \cite{wang2009fast,kogel2011fast,jerez2014embedded}. The system consists of a sequence of \(N\) masses connected by springs to each other, and to walls on either side. Each mass has value \(1\), and each spring has a spring constant of \(1\). See Fig.~\ref{fig: mass-spring model} for an illustration.

We model the dynamics of the oscillating masses system as follows. At time \(t\Delta\),
we let \(x_t=\begin{bmatrix}
    r_t^\top & s_t^\top
\end{bmatrix}^\top \) denote the state of the system, where the \(i\)-th element of vector \(r_t\in\mathbb{R}^N\) and \(s_t\in\mathbb{R}^N\) is the displacement and velocity of the \(i\)-th mass, respectively. Further, we let \(u_t\in\mathbb{R}^N\) denote the input to the system at time \(t\), whose \(i\)-th element is the external force exerted to the \(i\)-th mass. We let \(x_0 = \mathbf{0}_{2N}\) be the state of the system at time \(0\). Let \(L_N\in\mathbb{R}^{N\times N}\) is a symmetric tri-diagonal matrix whose diagonal entries are \(2\), and its sub-diagonal and super-diagonal entries are \(-1\). The discrete time dynamics of this system with sampling time period \(\Delta\) is given by \eqref{sys: DT ATI} and \eqref{eqn: ZOH} where 
\begin{equation*}
    A_c = \begin{bmatrix}
    0_{N\times N} &  I_N\\
     -L_N & 0_{N\times N}
    \end{bmatrix},\enskip B_c=\begin{bmatrix}
    0_{N\times N}\\
     I_N
    \end{bmatrix},\enskip h_c=\mathbf{0}_{2N}.
\end{equation*}

We consider the following constraints at each discrete time \(t\). The displacement, velocity and external force on each mass cannot exceed \([-\delta_1, \delta_1]\), \([-\delta_2, \delta_2]\) and \([-\rho, \rho]\), respectively, where \(\delta_1, \delta_2,\rho\in\mathbb{R}_+\). Further, for each external force, the maximum change in its magnitude within a sampling period \(\Delta\) is \(\gamma\). The aforementioned constraints are given by \eqref{mpc: input rate}, \eqref{mpc: state} and \eqref{mpc: input} where
\begin{equation}\label{eqn: mass param}
    \begin{aligned}
\mathbb{X}&= \{ r\in \mathbb{R}^N |  \norm{r}_\infty\leq \delta_1 \}\times\{ s\in \mathbb{R}^N | \norm{s}_\infty \leq \delta_2 \},\\
\mathbb{U} &= \{u\in\mathbb{R}^N| \norm{u}_\infty \leq \rho\}.
\end{aligned}
\end{equation}
Here the conic constraints in \eqref{mpc: state} and \eqref{mpc: input} (\ie, \(C_tx_t-a_t\geq 0\) and \(D_tx_t-b_t\geq 0\)) are not considered.

\subsubsection{Quadrotor path planning}

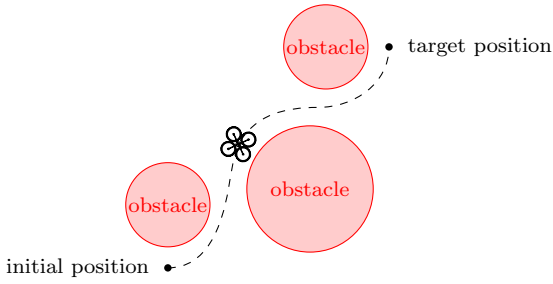
\begin{figure}[!ht]
	\centering
	\begin{tikzpicture}[scale=0.7]
		
		\coordinate (x0) at (-1.5, -2.7);
		\coordinate (xf) at (2.7, 1.5);
		\coordinate (o1) at (-1.5, -1.5);
		\coordinate (o2) at (1.2, -1.2);
		\coordinate (o3) at (1.5, 1.5);
		\coordinate (quad) at (-0.16,-0.35);
		\coordinate (mid1) at (-0.35, -1.2);
		\coordinate (mid2) at (1.35, 0.35);
		
		
		\node [quadcopter top,fill=white,draw=black,minimum width=1cm,rotate=-20,scale=0.4] at (quad) {}; 
		
		\draw[dashed] (x0) to[out=0,in=-100] (mid1) to[out=80,in=-110] (quad) to[out=60,in=180] (mid2) to[out=0,in=-90] (xf);
		
		\filldraw[color=red, fill=red!20](o1) circle (0.8) node {\scriptsize obstacle};
		\filldraw[color=red, fill=red!20](o2) circle (1.2) node {\scriptsize obstacle};
		\filldraw[color=red, fill=red!20](o3) circle (0.8) node {\scriptsize obstacle};
		
		\fill (x0) circle [radius=2pt];
		\node[label=left:{\scriptsize initial position}] at (x0) {};
		
		\fill (xf) circle [radius=2pt];
		\node[label=right:{\scriptsize target position}] at (xf) {};

	\end{tikzpicture}
		\caption{The quadrotor path planning problem.}
		\label{fig: quad path}
\end{figure}

We consider the problem of flying a quadrotor from its initial position to a target position while avoiding collision with cylinderical obstacles, see Fig.~\ref{fig: quad path} for an illustration. For the quadrotor dynamics, we consider the 3DoF model of the Autonomous Control Laboratory (ACL) custom quadrotor \cite[Ch.3]{szmuk2019successive}; see Fig.~\ref{fig: quad} and Fig.~\ref{fig: quad} for an illustration.

\begin{figure}[!ht]
    \centering
    \includegraphics[width=0.6\linewidth]{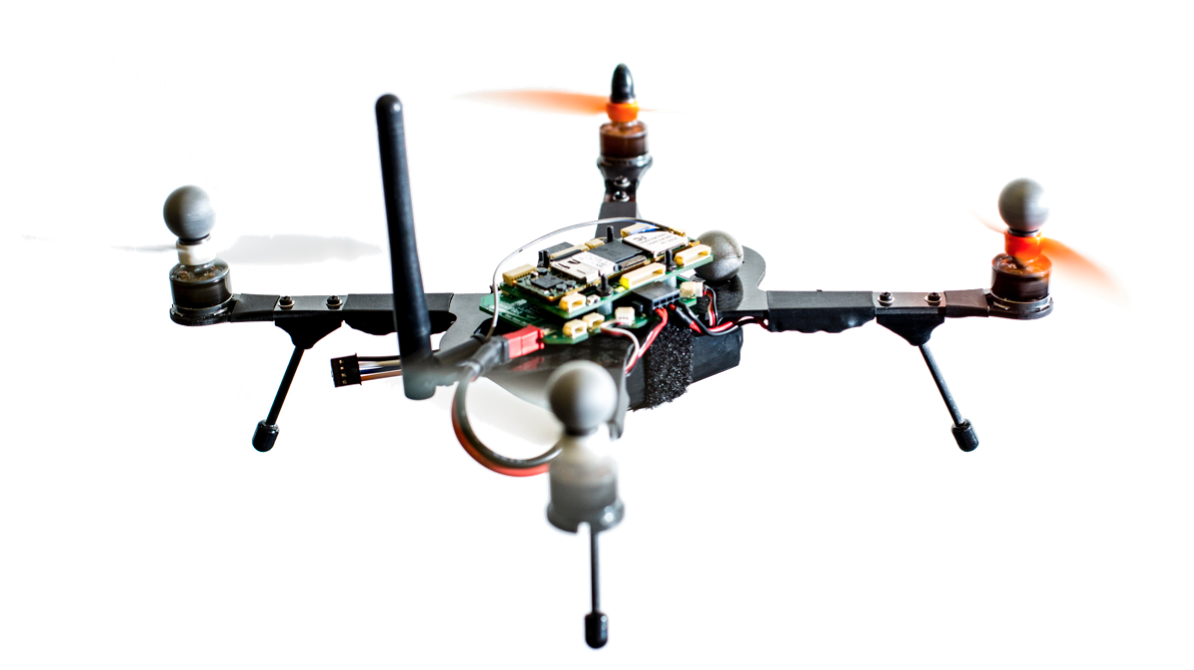}
    \caption{Autonomous Control Laboratory custom quadrotor}
    \label{fig: miki quad}
\end{figure}

We model the dynamics of the quadrotor as follows. At time \(t\Delta\), the state of the quadrotor is given by \(x_t=\begin{bmatrix}
r_t^\top & s_t^\top 
\end{bmatrix}^\top\), where \(r_t\in\mathbb{R}^3\) and \(s_t\in\mathbb{R}^3\) denote the position and velocity of the quadrotor's center of mass, respectively. We let \(x_0\) be the state of the system at time \(0\). The input of the quadrotor at time \(t\) is the thrust vector generated by its propellers, denoted by \(u_t\in\mathbb{R}^3\). Let \(m_0=0.35\) be the mass of the quadrotor and \(g_0=9.8\) be the gravitational constant. The discrete time quadrotor dynamics with sampling time period \(\Delta\) is given by \eqref{sys: DT ATI} and \eqref{eqn: ZOH} where
\begin{equation*}
    A_c = \begin{bmatrix}
    0_{3\times 3} &  I_3\\
    0_{3\times 3} & 0_{3\times 3}
    \end{bmatrix},\enskip B_c=\frac{1}{m_0}\begin{bmatrix}
    0_{3\times 3}\\
     I_3
    \end{bmatrix}, \enskip h_c=\begin{bmatrix}
    \mathbf{0}_5 \\ -g_0
    \end{bmatrix}.
\end{equation*}

We consider the following constraints on the thrust vector of the quadrotor. Due to the maximum power limit of the on-board motors, the magnitude of the thrust vector is upper bounded by \(\rho_1\in\mathbb{R}_+\). In addition, the vertical component of the thrust vector is lower bounded by \(\rho_2\in\mathbb{R}_+\) so the on-board motors are never turned off during the flight. The elementwise difference between two consecutive thrust vectors is upper bounded by \(\gamma\in\mathbb{R}_+\) to ensure a smooth thrust trajectory. Further, to upper bound the tilting angle of the quadrotor, we let the thrust vector be confined to a vertical icecream cone with half-angle \(\theta\in[0, \frac{\pi}{2}]\). Let \(e=\begin{bmatrix}0 & 0 & 1\end{bmatrix}^\top\), we can write the aforementioned constraints as \eqref{mpc: input rate} and \eqref{mpc: input} where
\begin{equation}\label{eqn: quad u consraints}
\begin{aligned}
    &\mathbb{U}= \left\{u\in\mathbb{R}^3\left| \norm{u}\cos\theta \leq\langle u, e \rangle, \norm{u}\leq \rho_1\right.\right\},\\
    &D_t=e, \enskip b_t=\rho_2.
\end{aligned}
\end{equation}

We also consider the following collision avoidance constraints. We want the position of the quadrotor to stay out of three vertical cylindrical region, \ie,
\begin{equation}\label{eqn: nonconvex 2}
  \norm{M x_t-o^i}^2\geq (\varrho^i)^2,
\end{equation}
for all \(i=1, 2, 3\), where \(M=\begin{bmatrix}
    I_2 & 0_{2\times 4}
\end{bmatrix}\), \(o^i\in\mathbb{R}^2\) and \(\varrho^i\in\mathbb{R}_+\) for all \(i=1, 2, 3\). However, the
above constraints are nonconvex, which render the resulting problem computationally challenging to solve. As a remedy, we consider the following linear approximation of \eqref{eqn: nonconvex 2}:
\begin{equation}\label{eqn: halfspace}
    \langle c_t^i, x_t\rangle \geq a_t^i,
\end{equation}
for all \(i=1, 2, 3\); see Fig.~\ref{fig: quad path} for an illustration. In Appendix~\ref{app: param}, we provide the detailed procedure on computing \(c_t^i\in\mathbb{R}^6\) and \(a_t^i\in\mathbb{R}\), and refer the interested reader to \cite[Sec. 4]{zagaris2018model} for a detailed discussion on this approximation.

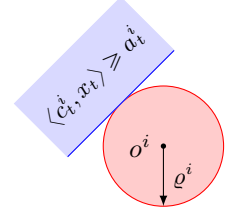
\begin{figure}[ht!]
\centering
  \begin{tikzpicture}[scale=0.4]

\filldraw[color=red, fill=red!20] (0,0) circle (2);

\filldraw[color=blue!15, fill=blue!15, rotate around={45:(0,0)}] (-2.5,2) rectangle (2.5, 4.5) node[pos=.5, rotate around={45:(0,0)}] {\color{black} \footnotesize $\langle c_t^i, x_t\rangle\geq a_t^i$};

\draw[color=blue, rotate around={45:(0,0)}] (-2.5, 2) -- (2.5, 2);

\draw[-latex] (0, 0) -- (0, -2) node[midway,right]{\footnotesize $\varrho^i$};

\node[circle,fill=black,inner sep=0pt,minimum size=2pt,label=left:{$o^i$}] (a) at (0,0) {};



\end{tikzpicture}
  \caption{Linearization of the constraint in \eqref{eqn: nonconvex 2} using \eqref{eqn: halfspace} }
  \label{fig: halfspace}
\end{figure}

With the above approximation, we can write the state constraints of the quadrotor in the form of \eqref{mpc: state} where 
\begin{equation}\label{eqn: quad x constraints}
\begin{aligned}
    &\mathbb{X}=\{r\in\mathbb{R}^3|  \norm{r}_\infty\leq \delta_1\}\times \{s\in\mathbb{R}^3| \norm{s}\leq \delta_2\}, \\
    & C_t=\begin{bmatrix}
    (c_t^1)^\top \\
    (c_t^2)^\top \\
    (c_t^3)^\top
    \end{bmatrix},\enskip  a_t=\begin{bmatrix}
    a_t^1\\
    a_t^2\\
    a_t^3
    \end{bmatrix}, \enskip \mathbb{K}_x=\mathbb{R}_+^3,
\end{aligned}
\end{equation}
for some \(\delta_1, \delta_2\in\mathbb{R}_+\). Here the set \(\mathbb{X}\) ensures the position and velocity of the quadrotor are bounded. 

\subsection{Numerical implementation and experiments}
\label{subsec: numerics}
We now discuss the numerical implementation of Algorithm~\ref{alg: PIPG} and demonstrate its performance using the two examples of constrained optimal control problems in Section~\ref{subsec: optimal control}.

\subsubsection{Efficient projections} The key step of implementing PIPG method is to compute the projection onto cone \(\mathbb{K}^\circ\) and set \(\mathbb{D}\). These projections can be computed efficiently for the following reasons. First, projections onto many common closed convex cones and sets can be computed using simple formulas, see \cite[Chp. 29]{bauschke2017convex} for some popular examples. The projection formula for the set \(\mathbb{U}\) in \eqref{eqn: quad u consraints} is given in \cite[Cor. 7.3]{bauschke2018projecting}. Second, let \(\mathbb{D}_1\subset\mathbb{R}^{n_1}\) and \(\mathbb{D}_2\subset\mathbb{R}^{n_2}\) be closed convex sets, \(x_1\in\mathbb{R}^{n_1}\) and \(x_2\in\mathbb{R}^{n_2}\). Then one can verify the following:
\begin{equation*}
    \pi_{\mathbb{D}_1\times \mathbb{D}_2}\left[\begin{bmatrix}
    x_1\\
    x_2
    \end{bmatrix}\right]=\begin{bmatrix}
    \pi_{\mathbb{D}_1}[x_1]\\
    \pi_{\mathbb{D}_2}[x_2]
    \end{bmatrix}.
\end{equation*}
Therefore, projections onto sets that are Cartesian products of sets with simple projection formulas, such as the set \(\mathbb{X}\) in \eqref{eqn: mass param} and \eqref{eqn: quad x constraints}, also admit simple formulas.

\subsubsection{Numerical experiments}

We demonstrate the numerical performance of PIPG using the two examples of optimization \eqref{opt: mpc}, namely the oscillating masses problem and the quadrotor path planning problem discussed in Section~\ref{subsec: optimal control}. We summarize the values of different problem parameters of these two examples in Appendix~\ref{app: param}. 

We compare the performance of PIPG, ADMM, PIPGeq and PDHG using optimization \eqref{opt: mpc} as follows. We initilize all methods using vectors whose entries are sampled from the standard normal distribution. We compare the performance of different methods using the convergence of the following two quantities:
\begin{equation}
    \text{error}_{\text{opt}}^j\coloneqq \frac{\norm{z^j-z^\star}^2}{\norm{z^\star}^2}, \enskip \text{error}_{\text{fea}}^j\coloneqq \frac{d_{\mathbb{K}}(Hz^j-g)}{\norm{z^\star}^2},
\end{equation}
where \(z^j\in\mathbb{D}\) is the candidate solution computed of optimization \eqref{opt: mpc} at the \(j\)-th iteration for \(j=2, 3, \ldots, k\), and \(z^\star\) be the ground truth optimal solution of optimization \eqref{opt: mpc} computed using commercial software Mosek \cite{mosek2019}. In addition, we also consider a restarting variant of PIPG where the iteration counter \(j\) is periodically reset to \(1\). Such restarting scheme is a popular heuristics for improving practical convergence performance of primal-dual methods \cite{su2016differential,xu2017accelerated}. 

\begin{figure}[!ht]
    \centering
    \begin{subfigure}[b]{0.49\columnwidth}
    \centering 
    \includegraphics[width=\textwidth]{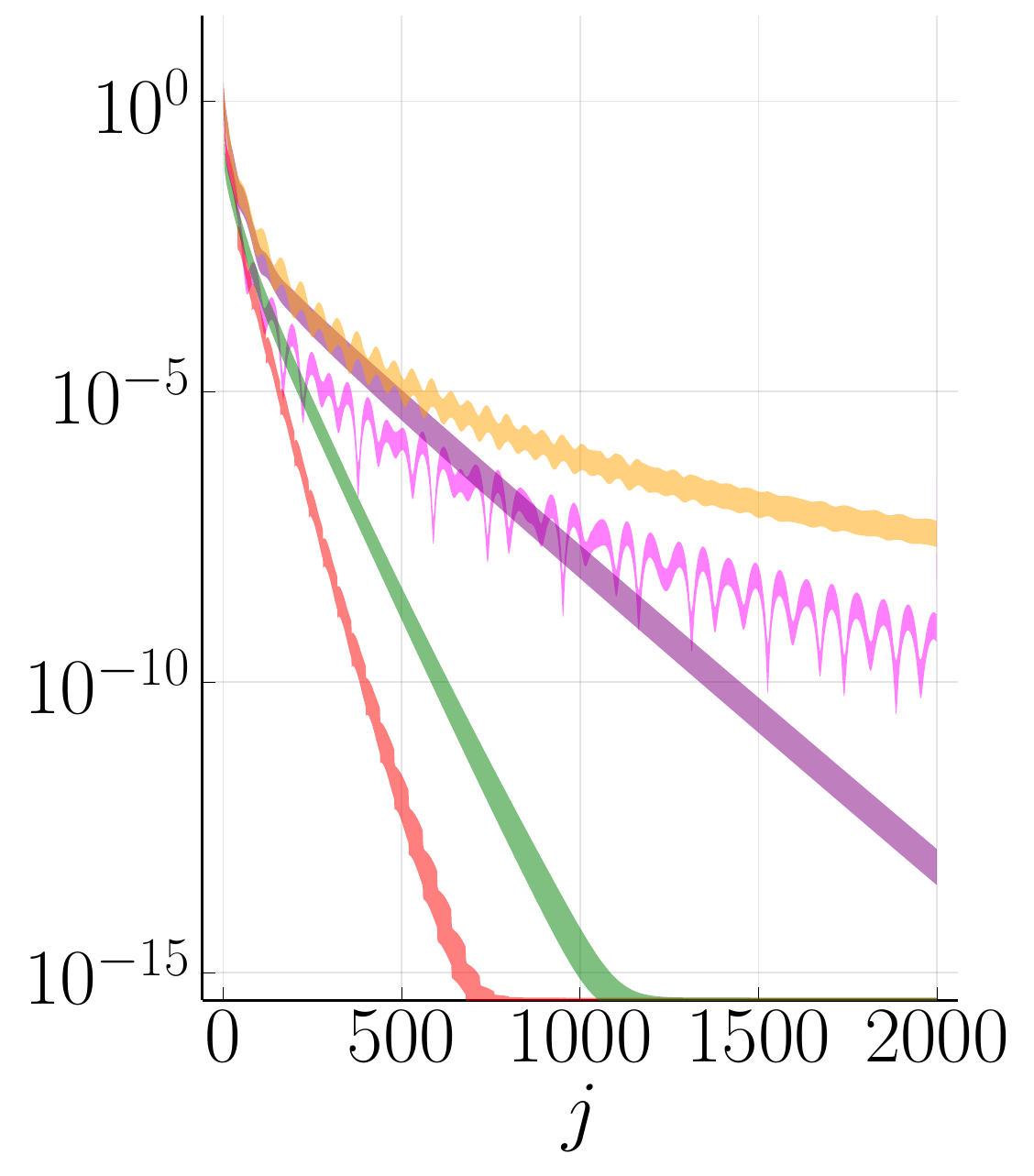}
   
    \caption{$\text{error}^j_{\text{opt}}$}
    \end{subfigure}
    \begin{subfigure}[b]{0.49\columnwidth}
    \centering 
    \includegraphics[width=\textwidth]{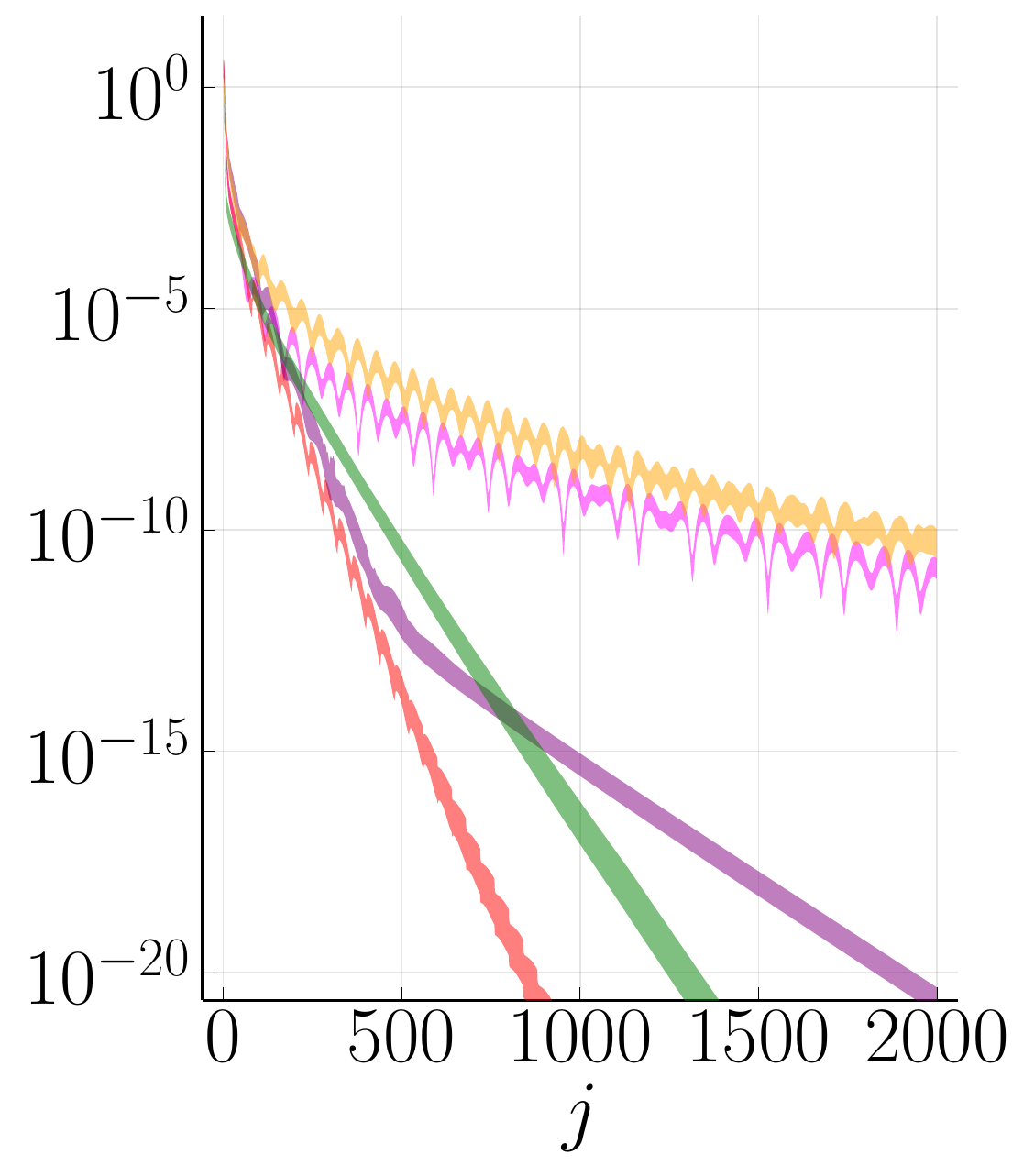}
  
    \caption{$\text{error}^j_{\text{fea}}$}
    \end{subfigure}

    \begin{subfigure}[b]{0.49\columnwidth}
    \centering 
    \includegraphics[width=\textwidth]{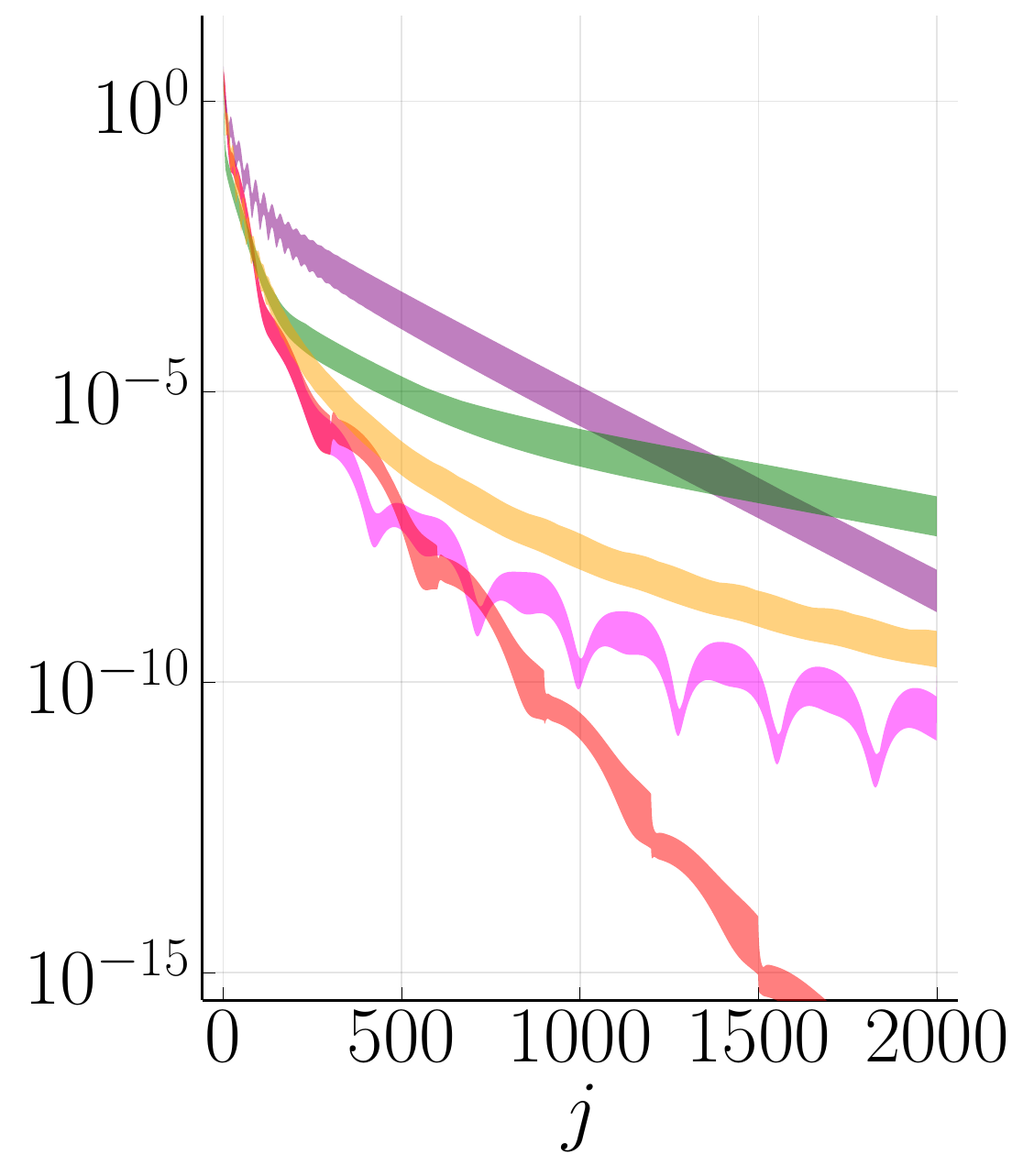}
  
    \caption{$\text{error}^j_{\text{opt}}$}
    \end{subfigure}
    \begin{subfigure}[b]{0.49\columnwidth}
    \centering 
    \includegraphics[width=\textwidth]{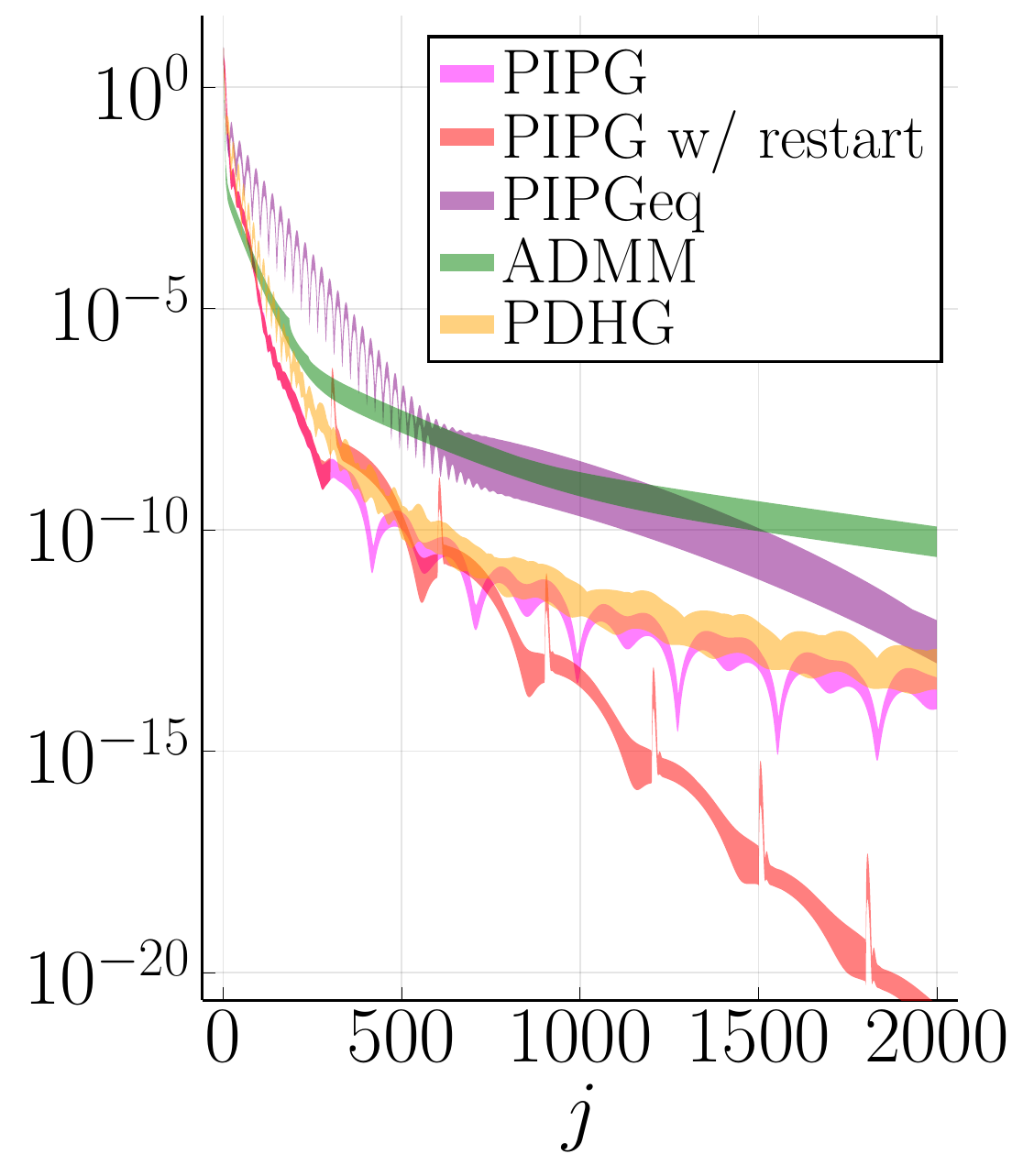}
  
    \caption{$\text{error}^j_{\text{fea}}$}
    \end{subfigure}
    \caption{Comparison of  different methods for oscillating masses problem (top row) and quadrotor path planning problem (bottom row). The shaded region shows the range of 100 different simulation results using independent random initializations. } 
    \label{fig: quad}
\end{figure}

The convergence results of different methods in terms of of \(e_{\text{opt}}^j\) and \(e_{\text{opt}}^j\) using 100 independent random initializations are illustrated in Fig.~\ref{fig: quad}. From these results we can see that PIPG clearly outperforms existing methods, especially when combined with the restarting heuristics. Note that, although the performance of ADMM is close to PIPG in the oscillating masses example, the per-iteration cost of ADMM is much higher than PIPG, as shown in Tab.~\ref{tab: comparison}. Therefore, PIPG still has clear advantage against ADMM.

\section{Conclusions}
\label{sec: conclusion}
We propose a novel primal-dual first-order method for conic optimization, named PIPG. We prove the convergence rates of PIPG in terms of the constraint violation and the primal-dual gap. We demonstrate the application of PIPG using examples in constrained optimal control. However, several questions still remain open. For example, it is unclear whether our method allow real-time implementation more efficient than interior point methods, or whether there are other restarting heuristcs better than the periodic one in Section~\ref{sec: numeric}. We aim to answer these open questions in our future work.




\appendix

\section{Proof of Proposition~\ref{prop: kkt}}
\label{app: prop}
We will use the following results.
\begin{lemma}\cite[Thm. 27.4]{rockafellar2015convex}\label{lem: set optimal} Let set \(\mathbb{D}\subset\mathbb{R}^n\) be closed and convex and function \(f:\mathbb{R}^n\to\mathbb{R}\) be continuously differentiable and convex. If  \(f(z^\star)\leq f(z)\) for all \(z\in\mathbb{D}\), then \(\langle \nabla f(z^\star), z-z^\star \rangle\geq 0\) for any \(z\in\mathbb{D}\).
\end{lemma}

\begin{lemma}\cite[Cor. 6.21]{rockafellar2009variational}\label{lem: polar cone}
If \(\mathbb{K}\subset\mathbb{R}^m\) is a closed convex cone, then \(\mathbb{K}^\circ\) is a closed convex cone and \((\mathbb{K}^\circ)^\circ=\mathbb{K}\).
\end{lemma}

We are now ready to prove Proposition~\ref{prop: kkt}.

\begin{pf*}{Proof}
First, if \eqref{eqn: suboptimality} holds, then we immediately have
\begin{equation}\label{eqn: suboptimal saddle}
    L(\overline{z}, w)\leq L(\overline{z}, \overline{w})\leq L(z, \overline{w})
\end{equation}
for all \(z\in\mathbb{D}\) and \(w\in \mathbb{K}^\circ\). The fist inequality above states that \(-L(\overline{z}, \overline{w})\leq -L(\overline{z}, w)\) for all \(w\in\mathbb{K}^\circ\), which, due to Lemma~\ref{lem: set optimal}, implies that
\begin{equation}\label{eqn: prop 1}
    \langle H\overline{z}-g, w-\overline{w}\rangle\leq 0
\end{equation}
for all \(w\in\mathbb{K}^\circ\). By letting \(w=0\) and \(w=2\overline{w}\) in \eqref{eqn: prop 1}, we conclude that
\begin{equation}\label{eqn: complementary}
    \langle H\overline{z}-g, \overline{w}\rangle=0.
\end{equation}
Combining \eqref{eqn: prop 1} and \eqref{eqn: complementary} gives \(\langle H\overline{z}-g, w\rangle\leq 0\) for all \(w\in\mathbb{K}^\circ\). Hence \(
    H\overline{z}-g\in(\mathbb{K}^\circ)^\circ=\mathbb{K}\), where the last step is due to Lemma~\ref{lem: polar cone}. 

Second, let \(z\) be such that \(z\in\mathbb{D}\) and \(Hz-g\in\mathbb{K}\). Since \(\overline{w}\in\mathbb{K}^\circ\), using \eqref{eqn: polar cone} we can show
\begin{equation}\label{eqn: prop 2}
    L(z, \overline{w})=f(z)+\langle Hz-g, \overline{w}\rangle\leq f(z).
\end{equation}
Further, using \eqref{eqn: suboptimal saddle} and \eqref{eqn: complementary} we can show
\begin{equation}\label{eqn: prop 3}
    f(\overline{z})=L(\overline{z}, \overline{w})\leq L(z, \overline{w}).
\end{equation}
By combining \eqref{eqn: prop 2} and \eqref{eqn: prop 3} we have \(f(\overline{z})\leq f(z)\). Since \(z\) is otherwise arbitrary except that \(z\in\mathbb{D}\) and \(Hz-g\in\mathbb{K}\), the proof is completed.
\end{pf*}
\section{Proof of Lemma~\ref{lem: iteration}}
\label{app: lem}
We start with some basic results that are necessary for the proof later. First, using \eqref{eqn: Bregman}, one can verify the following identity:
\begin{equation}\label{eqn: 3pt}
    \begin{aligned}
    &\langle \nabla f(z)-\nabla f(z'), z''-z\rangle\\
    &=B_f(z'', z')-B_f(z'',z)-B_f(z, z'), \enskip \forall z, z', z''\in\mathbb{R}^n.
    \end{aligned}
\end{equation}
If \(f=\norm{\cdot}^2\), the above identify becomes the following: 
\begin{equation}\label{eqn: los}
 2\langle z-z', z''-z\rangle=\norm{z''-z'}^2-\norm{z''-z}^2-\norm{z-z'}^2.
\end{equation}
Second, we will use Lemma~\ref{lem: polar cone}, together with the following existing results.

\begin{lemma}\cite[Lem. 2.2.7]{nesterov2018lectures}\label{lem: projection}
If set \(\mathbb{D}\subset\mathbb{R}^n\) is closed and convex, then \(\langle \pi_{\mathbb{D}}[z]-z, z'-\pi_{\mathbb{D}}[z]\rangle\geq 0\) for any \(z\in\mathbb{R}^n\) and \(z'\in\mathbb{D}\).
\end{lemma}

\begin{lemma}\cite[Thm. 6.30]{bauschke2017convex}\label{lem: Moreau}
If \(\mathbb{K}\subset\mathbb{R}^m\) is a closed convex cone, then \(\pi_{\mathbb{K}}[w]+\pi_{\mathbb{K}^\circ}[w]=w\) for all \(w\in\mathbb{R}^m\).
\end{lemma}

We are now ready to prove Lemma~\ref{lem: iteration}.
\begin{pf*}{Proof}
Let \(z, w, j\) be an arbitrary element in set \(\mathbb{D}\), cone \(\mathbb{K}^\circ\), and set \(\{1, 2, \ldots k\}\), respectively. We start with constructing an upper bound for \(L(z^{j+1}, w)-L(z, w^{j+1})\). To this end, first we use \eqref{eqn: Lagrangian} and \eqref{eqn: Bregman} to show the following identities
\begin{equation}\label{eqn: primal gap}
\begin{aligned}
 &L(z^{j+1}, w)-L(z, w)\\
 &=B_f(z^{j+1}, z)+\langle \nabla f(z)+H^\top w, z^{j+1}-z\rangle,
\end{aligned}
\end{equation}
\begin{equation}\label{eqn: dual gap}
L(z, w)-L(z, w^{j+1})=\langle Hz-g, w-w^{j+1}\rangle.
\end{equation}

Second, by applying Lemma~\ref{lem: projection} to the two projections in line~\ref{alg: w} and \ref{alg: z} in Algorithm~\ref{alg: PIPG} we can show the following two inequalities
\begin{subequations}
\begin{align}
    &0\leq \langle w^{j+1}-v^j-\beta^j(Hz^j-g), w-w^{j+1}\rangle, \label{eqn: proj v}\\
    &0\leq \langle z^{j+1}-z^j+\alpha^j(\nabla f(z^j)+H^\top w^{j+1}), z-z^{j+1} \rangle,\label{eqn: proj z}
\end{align}
\end{subequations}
Third, line~\ref{alg: v} in Algorithm~\ref{alg: PIPG} implies the following
\begin{equation}
    0= \langle v^{j+1}-w^{j+1}-\beta^jH(z^{j+1}-z^j), w-v^{j+1}\rangle.\label{eqn: proj w}
\end{equation}
Summing up \eqref{eqn: primal gap}, \eqref{eqn: dual gap}, \(\frac{1}{\beta^j}\times \)\eqref{eqn: proj v}, \(\frac{1}{\alpha^j}\times\)\eqref{eqn: proj z} and \(\frac{1}{\beta^j}\times\)\eqref{eqn: proj w} gives the following inequality
\begin{equation}\label{eqn: inner product}
    \begin{aligned}
        &L(z^{j+1}, w)-L(z, w^{j+1})\\ &\leq B_f(z^{j+1}, z)+\langle \nabla f(z)-\nabla f(z^j), z^{j+1}-z\rangle\\
        &\textstyle +\frac{1}{\alpha^j}\langle z^{j+1}-z^j, z-z^{j+1}\rangle\textstyle+\frac{1}{\beta^j}\langle w^{j+1}-v^j, w-w^{j+1}\rangle\\
        &\textstyle +\frac{1}{\beta^j}\langle v^{j+1}-w^{j+1}, w-v^{j+1}\rangle\\
        &+\langle v^{j+1}-w^{j+1}, H(z^{j+1}-z^j)\rangle.
    \end{aligned}
\end{equation}
Our next step is to bound the inner product terms in \eqref{eqn: inner product}. To this end, first we use \eqref{eqn: 3pt} and \eqref{eqn: los} to show the following identities
\begin{equation}\label{eqn: 3pt f}
    \begin{aligned}
        &\langle \nabla f(z)-\nabla f(z^j), z^{j+1}-z\rangle\\
        &= B_f(z^{j+1}, z^j)-B_f(z^{j+1}, z)-B_f(z, z^j),
    \end{aligned}
\end{equation}
\begin{equation}\label{eqn: 3pt z}
    \begin{aligned}
        &2\langle z^{j+1}-z^j, z-z^{j+1}\rangle\\
        &= \norm{z^j-z}^2-\norm{z^{j+1}-z}^2-\norm{z^{j+1}-z^j}^2, 
    \end{aligned}
\end{equation}
\begin{equation}\label{eqn: 3pt wv}
    \begin{aligned}
        &2\langle w^{j+1}-v^j, w-w^{j+1}\rangle\\
        &=\norm{v^j-w}^2-\norm{w^{j+1}-w}^2-\norm{w^{j+1}-v^j}^2,
    \end{aligned}
\end{equation}
\begin{equation}\label{eqn: 3pt w}
    \begin{aligned}
        &2\langle v^{j+1}-w^{j+1}, w-v^{j+1}\rangle\\
        &=\norm{w^{j+1}-w }^2-\norm{v^{j+1}-w}^2-\norm{v^{j+1}-w^{j+1}}^2.
    \end{aligned}
\end{equation}
Second, by completing the square we can show
\begin{equation}\label{eqn: square}
\begin{aligned}
&2\beta^j\langle v^{j+1}-w^{j+1}, H(z^{j+1}-z^j)\rangle\\
&\leq  \norm{v^{j+1}-w^{j+1}}^2+(\beta^j)^2\norm{H(z^{j+1}-z^j)}^2.
\end{aligned}
\end{equation}
Notice that now all inner product terms in \eqref{eqn: inner product} can be upper bounded. Finally, we further simplify these upper bounds. To this end, first we use the item 1 in Assumption~\ref{asp: opt} and the fact that \(\mnorm{H}^2\leq \sigma\) to show the following
\begin{subequations}
\begin{align}
    B_f(z^{j+1}, z^j)\leq \textstyle \frac{\lambda}{2}\norm{z^{j+1}-z^j}^2,\label{eqn: quad upper}\\
    -B_f(z, z^j)\leq \textstyle-\frac{\mu}{2}\norm{z^j-z}^2,\label{eqn: quad lower}\\
    \norm{H(z^{j+1}-z^j)}^2\leq \sigma \norm{z^{j+1}-z^j}^2.\label{eqn: matrix norm}
\end{align}
\end{subequations}
Second, we let \(y^j\coloneqq \frac{1}{\beta^j}(v^j+\beta^j(Hz^j-g)-w^{j+1})\). Applying Lemma~\ref{lem: polar cone} and Lemma~\ref{lem: Moreau} to the projection in line~\ref{alg: w} of Algorithm~\ref{alg: PIPG} we can show that \(\beta^j y^j\in(\mathbb{K}^\circ)^\circ=\mathbb{K}\). Since \(\mathbb{K}\) is a cone and \(\beta^j>0\), we know \(y^j\in\mathbb{K}\). Therefore, using \eqref{eqn: dist func} and definition of \(y^j\) we can show
\begin{equation}\label{eqn: dist2cone}
\begin{aligned}
     d_{\mathbb{K}}(Hz^j-g)\leq &\textstyle  \frac{1}{2}\norm{Hz^j-g-y^j}^2\\
     =& \textstyle\frac{1}{2(\beta^j)^2}\norm{w^{j+1}-v^j}^2.
\end{aligned}
\end{equation}
Finally, summing up \eqref{eqn: inner product}, \eqref{eqn: 3pt f}, \(\frac{1}{2\alpha^j}\times\)\eqref{eqn: 3pt z}, \(\frac{1}{2\beta^j}\times\)\eqref{eqn: 3pt wv}, \(\frac{1}{2\beta^j}\times\)\eqref{eqn: 3pt w}, 
\(\frac{1}{2\beta^j}\times\)\eqref{eqn: square},
\eqref{eqn: quad upper}, \eqref{eqn: quad lower}, \(\frac{\beta^j}{2}\times \)\eqref{eqn: matrix norm},  and \(\beta^j\times\)\eqref{eqn: dist2cone}, and using the assumption that \(\alpha^j(\lambda+\sigma\beta^j)=1\) we obtain the desired results.
\end{pf*}
\section{Proof of Theorem~\ref{thm: convex}}
\label{app: thm 1}
We will use the following result.
\begin{lemma}\cite[Lem. 3.1.1]{nesterov2018lectures}\label{lem: Jensen}
If function \(f:\mathbb{R}^n\to\mathbb{R}\) is convex, then 
\begin{equation}
   \textstyle f\left(\frac{1}{\sum_{j=1}^k\gamma^j}\sum_{j=1}^k\gamma^kz^j\right)\leq \frac{1}{\sum_{j=1}^k\gamma^j}\sum_{j=1}^k\gamma^jf(z^j)
\end{equation}
for any \(z^1, z^2, \ldots, z^k\in\mathbb{R}^n\) and \(\gamma^1, \gamma^2, \ldots, \gamma^k\in\mathbb{R}_+\). 
\end{lemma}

We are now ready to prove Theorem~\ref{thm: convex}.

\begin{pf*}{Proof}
Let \(z, w, j\) be an arbitrary element in set \(\mathbb{D}\), \(\mathbb{K}^\circ\) and \(\{1, 2, \ldots k\}\), respectively. Let \(V^j(z, w)= \frac{1}{2\alpha}\norm{z^j-z}^2+\frac{1}{2\beta}\norm{w^j-w}^2\). Since \(\alpha^j=\frac{1}{\beta\sigma+\lambda}\) and \(\beta^j=\beta\), the inequality in Lemma~\ref{lem: iteration} implies the following:
\begin{equation*}
\begin{aligned}
&L(z^{j+1}, w)-L(z, w^{j+1})+\beta d_{\mathbb{K}}(Hz^j-g)\\
&\leq V^j(z, w)-V^{j+1}(z, w),
\end{aligned}
\end{equation*}
for all \(z\in\mathbb{D}\), \(w\in \mathbb{K}^\circ\), and \(j=1, 2, \ldots, k\). Summing up this inequality for \(j=1, \ldots, k\) gives
\begin{equation}\label{eqn: thm1 telescope}
\begin{aligned}
&\textstyle \sum_{j=1}^k \big(L(z^{j+1}, w)-L(z, w^{j+1})+\beta d_{\mathbb{K}}(Hz^j-g)\big)\\
&\leq V^1(z, w)-V^{k+1}(z, w)\leq V^1(z, w),
\end{aligned}
\end{equation}
for all \(z\in\mathbb{D}\) and \(w\in \mathbb{K}^\circ\),
where the last step is because \(V^{k+1}(z, w)\geq 0\). From \eqref{eqn: dist func} and item 3 Assumption~\ref{asp: opt} we know that \(d_{\mathbb{K}}(Hz^j-g)\) and \(L(z^{j+1}, w^\star)-L(z^\star, w^{j+1})\) are non-negative for all \(j\). Hence \eqref{eqn: thm1 telescope} implies the following
\begin{equation*}
\begin{aligned}
\textstyle \sum_{j=1}^k \big(L(z^{j+1}, w)- L(z, w^{j+1})\big)\leq & V^1(z, w),\\
\textstyle \beta \sum_{j=1}^k d_{\mathbb{K}}(Hz^j-g)\leq & V^1(z^\star, w^\star),
\end{aligned}
\end{equation*}
for all \(z\in\mathbb{D}\), \(w\in \mathbb{K}^\circ\), where the second inequality is obtained by letting \(z=z^\star\) and \(w=w^\star\) in \eqref{eqn: thm1 telescope}.

Finally, applying the Jensen's inequality in \eqref{lem: Jensen} to convex function \(L(\cdot, w)\), \(-L(z, \cdot)\), and \(d_{\mathbb{K}}(\cdot)\) in the above two inequalities, respectively, we obtain the desired results.

\end{pf*}
\section{Proof of Theorem~\ref{thm: strongly convex}}
\label{app: thm 2}
We will use Lemma~\ref{lem: Jensen} in the following proof.
\begin{pf*}{Proof}
Let \(z, w, j\) be an arbitrary element in set \(\mathbb{D}\), \(\mathbb{K}^\circ\) and \(\{1, 2, \ldots k\}\), respectively. Let \(V^j(z, w)= \frac{1}{2\alpha^{j-1}}\norm{z^j-z}^2+\frac{1}{2\beta^{j-1}}\norm{v^j-w}^2\). Since \(\alpha^j=\frac{2}{(j+1)\mu+2\lambda}\) and \(\beta^j=\frac{(j+1)\mu}{2\sigma}\), the inequality in Lemma~\ref{lem: iteration} implies the following:
\begin{equation}\label{eqn: thm2 eqn1}
\begin{aligned}
    &\textstyle L(z^{j+1}, w)-L(z, w^{j+1})+\frac{(j+1)\mu}{2\sigma} d_{\mathbb{K}}(Hz^j-g)\\
    & \textstyle \leq  \frac{1}{2}(\frac{1}{\alpha^j}-\mu)\norm{z^j-z}^2+\frac{1}{2\beta^j}\norm{v^j-w}^2-V^{j+1}(z, w),
\end{aligned}    
\end{equation}
for all \(z\in\mathbb{D}\), \(w\in \mathbb{K}^\circ\), and \(j=1, 2, \ldots, k\). Let \(\kappa=\lambda/\mu\geq 1\), then one can verify the following 
\begin{equation}\label{eqn: thm2 eqn2}
\begin{aligned}
     \textstyle (\frac{1}{\alpha^j}-\mu)(j+2\kappa)= & \textstyle \frac{1}{\alpha^{j-1}}(j+2\kappa-1),\\
     \textstyle \frac{1}{\beta^j}(j+2\kappa)\leq & \textstyle \frac{1}{\beta^{j-1}}(j+2\kappa-1).
\end{aligned}
\end{equation}
Hence multiplying \eqref{eqn: thm2 eqn1} with \((j+2\kappa)\) then substituting in \eqref{eqn: thm2 eqn2} we can show
\begin{equation*}
\begin{aligned}
     &(j+2\kappa) (L(z^{j+1}, w)-L(z, w^{j+1}))\\
     &\textstyle + \frac{(j+1)(j+2\kappa)\mu}{2\sigma} d_{\mathbb{K}}(Hz^j-g)\\
     & \leq  (j+2\kappa-1)V^j(z, w)-(j+2\kappa)V^{j+1}(z, w),
\end{aligned}
\end{equation*}
for all \(z\in\mathbb{D}\), \(w\in \mathbb{K}^\circ\), and \(j=1, 2, \ldots, k\).
Summing up this inequality for \(j=1, 2, \ldots, k\) gives
\begin{equation}\label{eqn: thm2 telescope}
\begin{aligned}
 &\textstyle \sum_{j=1}^k (j+2\kappa) (L(z^{j+1}, w)-L(z, w^{j+1}))\\
 &+\textstyle \sum_{j=1}^k \frac{(j+1)(j+2\kappa)\mu}{2\sigma} d_{\mathbb{K}}(Hz^j-g)\\
 &\leq 2\kappa V^1(z, w)-(k+2\kappa)V^{k+1}(z, w) \leq 2\kappa V^1(z, w),
\end{aligned}
\end{equation}
for all \(z\in\mathbb{D}\) and \(w\in \mathbb{K}^\circ\), where the last step is because \(V^{k+1}(z, w)\geq 0\). From \eqref{eqn: dist func} and item 3 in Assumption~\ref{asp: opt} we know that \(d_{\mathbb{K}}(Hz^j-g)\) and \(\ell(z^{j+1}, w^{j+1}; z^\star, w^\star)\) are non-negative for all \(j\). Hence the above inequality implies the following
\begin{equation*}
    \begin{aligned}
     &\textstyle \sum_{j=1}^k (j+2)\big(L(z^{j+1}, w)- L(z, v^{j+1})\big)\leq 2\kappa V^1(z, w),\\
     &\textstyle \sum_{j=1}^k\frac{(j+1)(j+2)\mu}{2\sigma} d_{\mathbb{K}}(Hz^j-g)\leq  2\kappa V^1(z^\star, w^\star), 
    \end{aligned}
\end{equation*}
for all \(z\in\mathbb{D}\) and \(w\in \mathbb{K}^\circ\), where we used the fact that \(\kappa\geq 1\), and the second inequality is obtained by letting \(z=z^\star\) and \(w=w^\star\) in \eqref{eqn: thm2 telescope}.

Finally, applying the Jensen's inequality in Lemma~\ref{lem: Jensen} to convex function \(L(\cdot, w)\), \(-L(z, \cdot)\), and \(d_{\mathbb{K}}(\cdot)\) in the above two inequalities, respectively, we obtain the desired results.
\end{pf*}
\section{Transformation from an optimal control problem to a conic optimization}
\label{app: matrices}
We will use the following notation. We let \(\otimes\) denotes the Kronecker product, and \((\mathbb{D})^\tau\) denotes the Cartesian product of \(\tau\) copies of set \(\mathbb{D}\).

The optimization in \eqref{opt: mpc} is a special case of \eqref{opt: conic} by letting
\begin{equation*}
\begin{aligned}
    &z = \begin{bmatrix}
x_1^\top & x_2^\top & \cdots & x_\tau^\top & u_0^\top & u_1^\top & \cdots & u_{\tau-1}^\top
\end{bmatrix}^\top,\\
&f(z)=\frac{1}{2}z^\top P z+\langle p, z\rangle,\\
    & P=\begin{bmatrix}
    I_\tau\otimes Q & 0\\
    0 & I_\tau\otimes R
    \end{bmatrix},\\
    & p=\begin{bmatrix}
    \hat{x}_1^\top & \hat{x}_2^\top & \cdots & \hat{x}_\tau^\top & \hat{u}_1^\top & \hat{u}_2^\top & \cdots & \hat{u}_\tau^\top
    \end{bmatrix}^\top,\\
    & H=\begin{bmatrix}
    \overline{A} & \overline{B}\\
    0 & \overline{E}\\
     0 & -\overline{E}\\
    \overline{C} &  0\\
    0 & \overline{D}
    \end{bmatrix}, \enskip g =\begin{bmatrix}
    \overline{h}\\
    -\gamma \mathbf{1}\\
    -\gamma \mathbf{1}\\
    \overline{a}\\
    \overline{b}
    \end{bmatrix},\\
    &\mathbb{K}=0_{n_x\tau}\times \mathbb{R}_+^{2n_u(\tau-1)+\tau(p_x+p_u)}, \enskip \mathbb{D}=\overline{\mathbb{X}}\times \overline{\mathbb{U}},
\end{aligned}
\end{equation*}
where 
\begin{equation*}
\begin{aligned}
    &\overline{A}=I_{\tau n_x}-\begin{bmatrix}
    0 & 0\\
    I_{\tau-1}\otimes A & 0
    \end{bmatrix},\enskip  \overline{B}=-I_\tau\otimes B, \\
    &\overline{E}=\begin{bmatrix}
    I_{(\tau-1)n_u} & 0
    \end{bmatrix}-\begin{bmatrix}
    0 & I_{(\tau-1)n_u}
    \end{bmatrix},\\
    & \overline{C}=\begin{bmatrix}
    C_1 & 0 & \cdots &  0\\
    0 & C_2 & \cdots & 0\\
    0 & \ddots & \ddots & \vdots \\
    0 & 0 & \cdots & C_\tau
    \end{bmatrix},\enskip
     \overline{D}=\begin{bmatrix}
    D_1 & 0 & \cdots &  0\\
    0 & D_2 & \cdots & 0\\
    0 & \ddots & \ddots & \vdots \\
    0 & 0 & \cdots & D_\tau
    \end{bmatrix},\\
    & \overline{a}=\begin{bmatrix}
    a_1^\top & a_2^\top & \cdots & a_\tau^\top
    \end{bmatrix}^\top, \enskip \overline{b}=\begin{bmatrix}
    b_1^\top & b_2^\top & \cdots & b_\tau^\top
    \end{bmatrix}^\top,\\
    & \overline{h}=\begin{bmatrix}
    (Ax_0+h_0)^\top & h_1^\top & \cdots & h_{\tau-1}^\top
    \end{bmatrix}^\top,\\
    &\overline{\mathbb{X}}=(\mathbb{X})^\tau, \enskip \overline{\mathbb{U}}=(\mathbb{U})^\tau.
\end{aligned}
\end{equation*}

\section{Parameters of the optimal control problems in Section~\ref{subsec: numerics} }
\label{app: param}
\paragraph*{Oscillating masses}
We let the number of masses to be \(N=4\). 
In \eqref{opt: mpc}, we let \(\tau=30\), \(Q = I_{2N}\), \(R=I_N\), \(\hat{u}_t=\mathbf{0}_N\) and \(\hat{x}_{t+1}=\begin{bmatrix}
        \mathbf{1}_N^\top & \mathbf{0}_N^\top
\end{bmatrix}^\top \) for all \(t=0, 1, \ldots, \tau-1\). We also let \(\gamma=0.5\) in \eqref{mpc: input rate}, \(\Delta=0.25\) in \eqref{eqn: ZOH}, and \(\delta_1=\delta_2=\rho=2\) in \eqref{eqn: mass param}. 

\paragraph*{Quadrotor path planning}
In \eqref{opt: mpc}, we let \(\tau=30\), \(Q = \begin{bmatrix}
        I_3 & 0\\
        0 & 2.5I_3
\end{bmatrix}\), \(R=0.5I_3\), \(\hat{u}_t=\mathbf{0}_3\), and \(\hat{x}_{t+1}=\begin{bmatrix}
\hat{r}_{t+1}^\top & \hat{s}_{t+1}^\top
\end{bmatrix}^\top\) where \(\hat{s}_{t+1}=\mathbf{0}_3\) and 
\begin{equation*}
\hat{r}_{t+1} = \frac{t+1}{\tau}\begin{bmatrix} 2.5\\
1.5 \\
0 
\end{bmatrix} + \left(1-\frac{t+1}{\tau}\right)\begin{bmatrix} -1.5 \\
-2.5 \\
0 
\end{bmatrix}.  
\end{equation*}
for all \(t=0, 1, 2, \ldots, \tau-1\). We also let \(\gamma=3\) in \eqref{mpc: input rate}, \(\Delta=0.25\) in \eqref{eqn: ZOH}. We let \(\theta=\pi/4,\delta_1=3,\delta_2=5,\rho_1=5,\rho_2=2\) in \eqref{eqn: quad u consraints}. For all \(t=1, \ldots, \tau\), we let 
\begin{equation*}
    c_t^i=2M^\top (\tilde{r}_t-o_t^i),\enskip a_t^i=\norm{\tilde{r}_t}^2+(\varrho^i)^2-\norm{o^i}^2,
\end{equation*}
for \(i=1, 2, 3\) in \eqref{eqn: quad x constraints}, where
\begin{equation*}
\begin{aligned}
     &o^1=\begin{bmatrix}-1.5 \\ -1.5\end{bmatrix},\enskip  o^2=\begin{bmatrix}1.2 \\ -1.2\end{bmatrix},\enskip 
        o^3=\begin{bmatrix}1.5 \\ 1.5\end{bmatrix},\enskip
        \varrho^1=0.8, \\
    & \varrho^2=1.2,\enskip \varrho^3=0.8, \enskip M=\begin{bmatrix}
    I_2 & 0_{2\times 4}
    \end{bmatrix}, 
\end{aligned}
\end{equation*}
and \(\tilde{r}_t\) is computed as follows. If \(\norm{\hat{r}_t-o^i}\geq \varrho^i\) for all \(i=1, 2, 3\), then \(\tilde{r}_t=\hat{r}_t\). If there exists \(i\in\{1, 2,3\}\) such that \(\norm{\hat{r}_t-o^i}< \varrho^i\), then \(\tilde{r}_t=o^i+\frac{\varrho^i}{\norm{\hat{r}_t-o^i}}(\hat{r}_t-o^i)\). One can verify that \(\hat{r}_t\neq o^i\) for all \(t=1, 2, \ldots, \tau\) and \(i=1, 2, 3\) and there exists at most one \(i\in\{1, 2, 3\}\) such that \(\norm{\hat{r}_t-o^i}< \varrho^i\). Hence the \(\tilde{r}_t\) computed in the above manner is well defined and unique.


\bibliographystyle{apalike}      
\bibliography{reference}

\begin{thebibliography}{}

\bibitem[Andersen et~al., 2003]{andersen2003implementing}
Andersen, E.~D., Roos, C., and Terlaky, T. (2003).
\newblock On implementing a primal-dual interior-point method for conic
  quadratic optimization.
\newblock {\em Math. Program.}, 95(2):249--277.

\bibitem[Andersen et~al., 2011]{andersen2011interior}
Andersen, M., Dahl, J., Liu, Z., Vandenberghe, L., Sra, S., Nowozin, S., and
  Wright, S. (2011).
\newblock Interior-point methods for large-scale cone programming.
\newblock {\em Optim. Mach. Learn.}, 5583.

\bibitem[Bauschke et~al., 2018]{bauschke2018projecting}
Bauschke, H.~H., Bui, M.~N., and Wang, X. (2018).
\newblock Projecting onto the intersection of a cone and a sphere.
\newblock {\em SIAM J. Optim.}, 28(3):2158--2188.

\bibitem[Bauschke and Combettes, 2017]{bauschke2017convex}
Bauschke, H.~H. and Combettes, P.~L. (2017).
\newblock {\em Convex Analysis and Monotone Operator Theory in Hilbert Spaces},
  volume 408.
\newblock Springer.

\bibitem[Ben-Tal and Nemirovski, 2001]{ben2001lectures}
Ben-Tal, A. and Nemirovski, A. (2001).
\newblock {\em Lectures on Modern Convex Optimization: Analysis, Algorithms,
  and Engineering Applications}.
\newblock SIAM.

\bibitem[Betts, 2010]{betts2010practical}
Betts, J. (2010).
\newblock {\em Practical Methods for Optimal Control and Estimation using
  Nonlinear Programming}.
\newblock SIAM, Philadelphia.

\bibitem[Boyd et~al., 2011]{boyd2011distributed}
Boyd, S., Parikh, N., and Chu, E. (2011).
\newblock {\em Distributed Optimization and Statistical Learning via the
  Alternating Direction Method of Multipliers}.
\newblock Now Publishers Inc.

\bibitem[Boyd and Vandenberghe, 2004]{boyd2004convex}
Boyd, S.~P. and Vandenberghe, L. (2004).
\newblock {\em Convex Optimization}.
\newblock Cambridge University Press.

\bibitem[Chambolle and Pock, 2011]{chambolle2011first}
Chambolle, A. and Pock, T. (2011).
\newblock A first-order primal-dual algorithm for convex problems with
  applications to imaging.
\newblock {\em J. Math. Imaging Vis.}, 40(1):120--145.

\bibitem[Chambolle and Pock, 2016a]{chambolle2016introduction}
Chambolle, A. and Pock, T. (2016a).
\newblock An introduction to continuous optimization for imaging.
\newblock {\em Acta Numerica}, 25:161--319.

\bibitem[Chambolle and Pock, 2016b]{chambolle2016ergodic}
Chambolle, A. and Pock, T. (2016b).
\newblock On the ergodic convergence rates of a first-order primal--dual
  algorithm.
\newblock {\em Math. Program.}, 159(1-2):253--287.

\bibitem[Chen, 1999]{chen1999linear}
Chen, C.-T. (1999).
\newblock {\em {L}inear {S}ystem {T}heory and {D}esign}.
\newblock {O}xford {U}niversity {P}ress, New York.

\bibitem[Chen et~al., 2013]{chen2013primal}
Chen, P., Huang, J., and Zhang, X. (2013).
\newblock A primal--dual fixed point algorithm for convex separable
  minimization with applications to image restoration.
\newblock {\em Inverse Problems}, 29(2):025011.

\bibitem[Chen et~al., 2016]{chen2016primal}
Chen, P., Huang, J., and Zhang, X. (2016).
\newblock A primal-dual fixed point algorithm for minimization of the sum of
  three convex separable functions.
\newblock {\em Fixed Point Theory Appl.}, 2016(1):1--18.

\bibitem[Cohen et~al., 2018]{cohen2018acceleration}
Cohen, M., Diakonikolas, J., and Orecchia, L. (2018).
\newblock On acceleration with noise-corrupted gradients.
\newblock In {\em Int. Conf. Mach. Learn.}, pages 1019--1028. PMLR.

\bibitem[Condat, 2013]{condat2013primal}
Condat, L. (2013).
\newblock A primal--dual splitting method for convex optimization involving
  lipschitzian, proximable and linear composite terms.
\newblock {\em J. Optim. Theory Appl.}, 158(2):460--479.

\bibitem[Davis and Yin, 2017]{davis2017three}
Davis, D. and Yin, W. (2017).
\newblock A three-operator splitting scheme and its optimization applications.
\newblock {\em Set-Valued Var. Anal.}, 25(4):829--858.

\bibitem[Eckstein, 1989]{eckstein1989splitting}
Eckstein, J. (1989).
\newblock {\em Splitting Methods for Monotone Operators with Applications to
  Parallel Optimization}.
\newblock PhD thesis, Massachusetts Inst. Technol.

\bibitem[Eren et~al., 2017]{eren2017model}
Eren, U., Prach, A., Ko{\c{c}}er, B.~B., Rakovi{\'c}, S.~V., Kayacan, E., and
  A{\c{c}}{\i}kme{\c{s}}e, B. (2017).
\newblock Model predictive control in aerospace systems: Current state and
  opportunities.
\newblock {\em J. Guid. Control Dyn.}, 40(7):1541--1566.

\bibitem[Fortin and Glowinski, 2000]{fortin2000augmented}
Fortin, M. and Glowinski, R. (2000).
\newblock {\em Augmented Lagrangian methods: Applications to the Numerical
  Solution of Boundary-Value Problems}.
\newblock Elsevier.

\bibitem[Gabay and Mercier, 1976]{gabay1976dual}
Gabay, D. and Mercier, B. (1976).
\newblock A dual algorithm for the solution of nonlinear variational problems
  via finite element approximation.
\newblock {\em Comput. Math. Appl.}, 2(1):17--40.

\bibitem[Goldstein et~al., 2014]{goldstein2014fast}
Goldstein, T., O'{D}onoghue, B., Setzer, S., and Baraniuk, R. (2014).
\newblock Fast alternating direction optimization methods.
\newblock {\em SIAM J. Imag. Sci.}, 7(3):1588--1623.

\bibitem[He and Yuan, 2012]{he20121}
He, B. and Yuan, X. (2012).
\newblock On the ${O}(1/n)$ convergence rate of the douglas--rachford
  alternating direction method.
\newblock {\em SIAM J Numer. Anal.}, 50(2):700--709.

\bibitem[Jerez et~al., 2014]{jerez2014embedded}
Jerez, J.~L., Goulart, P.~J., Richter, S., Constantinides, G.~A., Kerrigan,
  E.~C., and Morari, M. (2014).
\newblock Embedded online optimization for model predictive control at
  megahertz rates.
\newblock {\em IEEE Trans. Automat. Control}, 59(12):3238--3251.

\bibitem[Kadkhodaie et~al., 2015]{kadkhodaie2015accelerated}
Kadkhodaie, M., Christakopoulou, K., Sanjabi, M., and Banerjee, A. (2015).
\newblock Accelerated alternating direction method of multipliers.
\newblock In {\em Proc Int. Conf. Knowl. Discovery Data Mining}, pages
  497--506.

\bibitem[K{\"o}gel and Findeisen, 2011]{kogel2011fast}
K{\"o}gel, M. and Findeisen, R. (2011).
\newblock Fast predictive control of linear systems combining {N}esterov's
  gradient method and the method of multipliers.
\newblock In {\em Proc. IEEE Conf. Decision Control and Eur. Control Conf.},
  pages 501--506. IEEE.

\bibitem[Korpelevich, 1977]{korpelevich1977extragradient}
Korpelevich, G. (1977).
\newblock Extragradient method for finding saddle points and other problems.
\newblock {\em Matekon}, 13(4):35--49.

\bibitem[Krol et~al., 2012]{krol2012preconditioned}
Krol, A., Li, S., Shen, L., and Xu, Y. (2012).
\newblock Preconditioned alternating projection algorithms for maximum a
  posteriori ect reconstruction.
\newblock {\em Inverse problems}, 28(11):115005.

\bibitem[Lan et~al., 2011]{lan2011primal}
Lan, G., Lu, Z., and Monteiro, R.~D. (2011).
\newblock Primal-dual first-order methods with \(\mathcal{O}(1/\epsilon)\)
  iteration-complexity for cone programming.
\newblock {\em Math. Program.}, 126(1):1--29.

\bibitem[Liu et~al., 2017]{liu2017survey}
Liu, X., Lu, P., and Pan, B. (2017).
\newblock Survey of convex optimization for aerospace applications.
\newblock {\em Astrodynamics}, 1(1):23--40.

\bibitem[Luo and Yu, 2006]{luo2006introduction}
Luo, Z.-Q. and Yu, W. (2006).
\newblock An introduction to convex optimization for communications and signal
  processing.
\newblock {\em IEEE J. Sel. Areas Commun.}, 24(8):1426--1438.

\bibitem[Majumdar et~al., 2020]{majumdar2020recent}
Majumdar, A., Hall, G., and Ahmadi, A.~A. (2020).
\newblock Recent scalability improvements for semidefinite programming with
  applications in machine learning, control, and robotics.
\newblock {\em Annu. Rev Control Robot. Auton. Syst.}, 3:331--360.

\bibitem[Malyuta et~al., 2021]{malyuta2021advances}
Malyuta, D., Yu, Y., Elango, P., and A{\c{c}}ikme{\c{s}}e, B. (2021).
\newblock Advances in trajectory optimization for space vehicle control.
\newblock {\em arXiv preprint arXiv:2108.02335 [math.OC]}.

\bibitem[{MOSEK ApS}, 2019]{mosek2019}
{MOSEK ApS} (2019).
\newblock {\em The MOSEK optimization toolbox for MATLAB manual. Version 9.0.}

\bibitem[Nemirovski, 2004]{nemirovski2004prox}
Nemirovski, A. (2004).
\newblock Prox-method with rate of convergence ${O}(1/t)$ for variational
  inequalities with {L}ipschitz continuous monotone operators and smooth
  convex-concave saddle point problems.
\newblock {\em SIAM J. Optim}, 15(1):229--251.

\bibitem[Nesterov, 2007]{nesterov2007dual}
Nesterov, Y. (2007).
\newblock Dual extrapolation and its applications to solving variational
  inequalities and related problems.
\newblock {\em Math. Program.}, 109(2):319--344.

\bibitem[Nesterov, 2018]{nesterov2018lectures}
Nesterov, Y. (2018).
\newblock {\em Lectures on Convex Optimization}, volume 137.
\newblock Springer.

\bibitem[Nesterov and Nemirovskii, 1994]{nesterov1994interior}
Nesterov, Y. and Nemirovskii, A. (1994).
\newblock {\em Interior-Point Polynomial Algorithms in Convex Programming}.
\newblock SIAM.

\bibitem[Ouyang et~al., 2015]{ouyang2015accelerated}
Ouyang, Y., Chen, Y., Lan, G., and Pasiliao~Jr, E. (2015).
\newblock An accelerated linearized alternating direction method of
  multipliers.
\newblock {\em SIAM J. Imag. Sci.}, 8(1):644--681.

\bibitem[O’Connor and Vandenberghe, 2020]{o2020equivalence}
O’Connor, D. and Vandenberghe, L. (2020).
\newblock On the equivalence of the primal-dual hybrid gradient method and
  douglas--rachford splitting.
\newblock {\em Math. Program.}, 179(1):85--108.

\bibitem[O’Donoghue et~al., 2016]{o2016conic}
O’Donoghue, B., Chu, E., Parikh, N., and Boyd, S. (2016).
\newblock Conic optimization via operator splitting and homogeneous self-dual
  embedding.
\newblock {\em J. Optim. Theory Appl.}, 169(3):1042--1068.

\bibitem[Rockafellar, 2015]{rockafellar2015convex}
Rockafellar, R.~T. (2015).
\newblock {\em Convex Analysis}.
\newblock Princeton University Press.

\bibitem[Rockafellar and Wets, 2009]{rockafellar2009variational}
Rockafellar, R.~T. and Wets, R. J.-B. (2009).
\newblock {\em Variational Analysis}, volume 317.
\newblock Springer Science \& Business Media.

\bibitem[Stellato et~al., 2020]{stellato2020osqp}
Stellato, B., Banjac, G., Goulart, P., Bemporad, A., and Boyd, S. (2020).
\newblock {OSQP}: an operator splitting solver for quadratic programs.
\newblock {\em Math. Program. Comput.}, 12(4):637--672.

\bibitem[Su et~al., 2016]{su2016differential}
Su, W., Boyd, S., and Candes, E.~J. (2016).
\newblock A differential equation for modeling nesterov's accelerated gradient
  method: Theory and insights.
\newblock {\em J. Mach. Learn. Res.}, 17(1):5312--5354.

\bibitem[Szmuk, 2019]{szmuk2019successive}
Szmuk, M. (2019).
\newblock {\em Successive Convexification \& High Performance Feedback Control
  for Agile Flight}.
\newblock PhD thesis, Dept. Aeronatu. \& Astronaut., Univ. Washington.

\bibitem[V{\~u}, 2013]{vu2013splitting}
V{\~u}, B.~C. (2013).
\newblock A splitting algorithm for dual monotone inclusions involving
  cocoercive operators.
\newblock {\em Adv. Comput. Math.}, 38(3):667--681.

\bibitem[Wang and Banerjee, 2014]{wang2014bregman}
Wang, H. and Banerjee, A. (2014).
\newblock Bregman alternating direction method of multipliers.
\newblock {\em Proc. Adv. Neural Inf. Process. Syst.}, 4(January):2816--2824.

\bibitem[Wang and Elia, 2010]{wang2010control}
Wang, J. and Elia, N. (2010).
\newblock Control approach to distributed optimization.
\newblock In {\em Proc. Allerton Conf. Commun. Control Comput.}, pages
  557--561. IEEE.

\bibitem[Wang and Boyd, 2009]{wang2009fast}
Wang, Y. and Boyd, S. (2009).
\newblock Fast model predictive control using online optimization.
\newblock {\em IEEE Trans. Control Syst. Technol.}, 18(2):267--278.

\bibitem[Xu, 2017]{xu2017accelerated}
Xu, Y. (2017).
\newblock Accelerated first-order primal-dual proximal methods for linearly
  constrained composite convex programming.
\newblock {\em SIAM J. Optim.}, 27(3):1459--1484.

\bibitem[Yan, 2018]{yan2018new}
Yan, M. (2018).
\newblock A new primal--dual algorithm for minimizing the sum of three
  functions with a linear operator.
\newblock {\em J. Sci. Comput.}, 76(3):1698--1717.

\bibitem[Yu and A{\c{c}}{\i}kme{\c{s}}e, 2020]{yu2020rlc}
Yu, Y. and A{\c{c}}{\i}kme{\c{s}}e, B. (2020).
\newblock {RLC} circuits-based distributed mirror descent method.
\newblock {\em IEEE Control Syst. Lett.}, 4(3):548--553.

\bibitem[Yu et~al., 2020a]{yu2020mass}
Yu, Y., A{\c{c}}{\i}kme{\c{s}}e, B., and Mesbahi, M. (2020a).
\newblock Mass--spring--damper networks for distributed optimization in
  non-{E}uclidean spaces.
\newblock {\em Automatica}, 112:108703.

\bibitem[Yu et~al., 2020b]{yu2020proportional}
Yu, Y., Elango, P., and A{\c{c}}{\i}kme{\c{s}}e, B. (2020b).
\newblock Proportional-integral projected gradient method for model predictive
  control.
\newblock {\em IEEE Control Syst. Lett.}

\bibitem[Zagaris et~al., 2018]{zagaris2018model}
Zagaris, C., Park, H., Virgili-Llop, J., Zappulla, R., Romano, M., and
  Kolmanovsky, I. (2018).
\newblock Model predictive control of spacecraft relative motion with
  convexified keep-out-zone constraints.
\newblock {\em J. Guid. Control Dyn.}, 41(9):2054--2062.

\end{thebibliography}

\end{document}